%%%%%%%%%%%%%%%%%%%%%%%%%%%%%%%%%%%%%%%%%%%%%%%%%%%%%%%%%%
%%%%%%%%%%%%%%%%%%%%%%%%%%%%%%%%%%%%%%%%%%%%%%%%%%%%%%%%%%
%%
%%     This is the AMS-LaTeX file:
%%
%%     Auricchio-Colli-Gilardi-Reali-Rocca 1
%%     modello sul covid
%%
%%%%%%%%%%%%%%%%%%%%%%%%%%%%%%%%%%%%%%%%%%%%%%%%%%%%%%%%%%

\let\Begin\begin
\let\End\end

%%%%%%%%%%%%%%%%%%%%%%%%%%%%%%%%%
%% page layout
%%%%%%%%%%%%%%%%%%%%%%%%%%%%%%%%%

\documentclass[twoside,12pt]{article}
\setlength{\textheight}{24cm}
\setlength{\textwidth}{16cm}
\setlength{\oddsidemargin}{2mm}
\setlength{\evensidemargin}{2mm}
\setlength{\topmargin}{-15mm}
\parskip2mm

%%%%%%%%%%%%%%%%%%%%%%%%%%%%%%%%%
%% packages
%%%%%%%%%%%%%%%%%%%%%%%%%%%%%%%%%

%\usepackage{color}
\usepackage[usenames,dvipsnames]{color}
\usepackage{amsmath}
\usepackage{amsthm}
\usepackage{amssymb}
\usepackage[mathcal]{euscript}

\usepackage[notref,notcite]{showkeys}
\usepackage{showkeys}
%
%		COLORS FOR CORRECTIONS
%
% do the same, please (i.e., don't use the standard {\color{red} text} or similar): 
% just choose the color you prefer in \def\yourname

% example of use:  \gianni{I want this to become blue}
% warning: no blank lines (\par) in the argument of \yourname

\definecolor{viola}{rgb}{0.3,0,0.7}
\definecolor{ciclamino}{rgb}{0.5,0,0.5}

\def\gianni #1{{\color{red}#1}}
\def\pier #1{{\color{blue}#1}}
\def\ale #1{{\color{magenta}#1}}
\def\pcol #1{{\color{red}#1}}

\def\oldgianni #1{#1}
\def\gianni #1{#1}
\def\pier #1{#1}
\def\ale #1{#1}
\def\pcol #1{#1}
%\def\revis #1{#1}

%%%%%%%%%%%%%%%%%%%%%%%%%%%%%%%%%
%% bibliographystyle
%%%%%%%%%%%%%%%%%%%%%%%%%%%%%%%%%

\bibliographystyle{plain}

%%%%%%%%%%%%%%%%%%%%%%%%%%%%%%%%%
%% environments
%%%%%%%%%%%%%%%%%%%%%%%%%%%%%%%%%

%

\def\Beq{\Begin{equation}}
\def\Eeq{\End{equation}}
\def\Bsist{\Begin{eqnarray}}
\def\Esist{\End{eqnarray}}

\def\Bthm{\Begin{theorem}}
\def\Ethm{\End{theorem}}
\def\Blem{\Begin{lemma}}
\def\Elem{\End{lemma}}
\def\Bprop{\Begin{proposition}}
\def\Eprop{\End{proposition}}

\def\Brem{\Begin{remark}\rm}
\def\Erem{\End{remark}}

\def\Bnot{\Begin{notation}\rm}
\def\Enot{\End{notation}}
\def\Bdim{\Begin{proof}}
\def\Edim{\End{proof}}
\def\Bcenter{\Begin{center}}
\def\Ecenter{\End{center}}
\let\non\nonumber

%%%%%%%%%%%%%%%%%%%%%%%%%%%%%%%%%
%% macros
%%%%%%%%%%%%%%%%%%%%%%%%%%%%%%%%%

% macro salvate

% sottosezioni non numerate

\def\step #1 \par{\medskip\noindent{\bf #1.}\quad}

% abbreviazioni di parole

\def\Lip{Lip\-schitz}
\def\Holder{H\"older}

\def\aand{\quad\hbox{and}\quad}

\def\lhs{left-hand side}
\def\rhs{right-hand side}
\def\sfw{straightforward}

\let\hat\widehat
\let\tilde\widetilde

% bold, cal e mathop

\def\multibold #1{\def\arg{#1}%
  \ifx\arg\pto \let\next\relax
  \else
  \def\next{\expandafter
    \def\csname #1#1#1\endcsname{{\bf #1}}%
    \multibold}%
  \fi \next}

\def\pto{.}

\def\multical #1{\def\arg{#1}%
  \ifx\arg\pto \let\next\relax
  \else
  \def\next{\expandafter
    \def\csname cal#1\endcsname{{\cal #1}}%
    \multical}%
  \fi \next}

% operatori

\def\multimathop #1 {\def\arg{#1}%
  \ifx\arg\pto \let\next\relax
  \else
  \def\next{\expandafter
    \def\csname #1\endcsname{\mathop{\rm #1}\nolimits}%
    \multimathop}%
  \fi \next}

\multibold
qwertyuiopasdfghjklzxcvbnmQWERTYUIOPASDFGHJKLZXCVBNM.

\multical
QWERTYUIOPASDFGHJKLZXCVBNM.

\multimathop
diag dist div dom mean meas sign supp .

% accorpamenti di formule citate:
% uso  \accorpa {prima}{seconda}
%      \Accorpa\cs prima seconda (con il comodo blank anche dopo)
% NB: \Accorpa definisce \cs come l'accorpamento delle due citazioni
% e scrive sul file.log

\def\accorpa #1#2{\eqref{#1}--\eqref{#2}}
\def\Accorpa #1#2 #3 {\gdef #1{\eqref{#2}--\eqref{#3}}%
  \wlog{}\wlog{\string #1 -> #2 - #3}\wlog{}}

% macro comode

\def\separa{\noalign{\allowbreak}}

\def\somma #1#2#3{\sum_{#1=#2}^{#3}}

\def\graffe #1{\mathopen\{#1\mathclose\}}

\def\<#1>{\mathopen\langle #1\mathclose\rangle}
\def\norma #1{\mathopen \| #1\mathclose \|}

\def\[#1]{\mathopen\langle\!\langle #1\mathclose\rangle\!\rangle}

\def\iot {\int_0^t}
\def\ioT {\int_0^T}
\def\intQt{\int_{Q_t}}
\def\intQ{\int_Q}
\def\iO{\int_\Omega}

\def\dt{\partial_t}
\def\dn{\partial_\nu}

\def\cpto{\,\cdot\,}

\def\checkmmode #1{\relax\ifmmode\hbox{#1}\else{#1}\fi}
\def\aeO{\checkmmode{a.e.\ in~$\Omega$}}
\def\aeQ{\checkmmode{a.e.\ in~$Q$}}
\def\aet{\checkmmode{a.e.\ in~$(0,T)$}}

% insiemi numerici

\def\erre{{\mathbb{R}}}

% spazi di funzioni a valori vettoriali su [0,T], [0,t], [0,s], [0,+\infty), [\delta,T]

% Come ricordare: in generale i simboli L H W  C da soli per gli spazi su (0,T)
% gli stessi raddoppiati per (0,+\infty)
% aggiunta di t o s al simbolo per (0,t) e (0,s)
% aggiunta di d al simbolo semplice o doppio per intervalli (\delta,T) e (\delta,+\infty)
% il simbolo C e i suoi derivati mettono le quadre anziche' le tonde

% Esempi   \L2V   \L\infty\Vp   \W{1,1}H   \C0H   \LL2V   \CC0\Vp   \Ld2V  \CCdH

\def\genspazio #1#2#3#4#5{#1^{#2}(#5,#4;#3)}
\def\spazio #1#2#3{\genspazio {#1}{#2}{#3}T0}

\def\L {\spazio L}
\def\H {\spazio H}
\def\W {\spazio W}

\def\C #1#2{C^{#1}([0,T];#2)}

% spazi di funzioni su \Omega, \Gamma, Q e \Sigma

\def\Lx #1{L^{#1}(\Omega)}
\def\Hx #1{H^{#1}(\Omega)}
\def\Wx #1{W^{#1}(\Omega)}

\def\Cx #1{C^{#1}(\overline\Omega)}

\def\LQ #1{L^{#1}(Q)}

\def\CQ #1{C^{#1}(\overline Q)}

\def\Ldue{\Lx 2}
\def\Linfty{\Lx\infty}

\def\Huno{\Hx 1}
\def\Hdue{\Hx 2}

% spazi di funzioni su Q e S

\def\LQ #1{L^{#1}(Q)}

% lettere greche

\let\theta\vartheta
\let\eps\varepsilon
\let\ka\kappa

\let\TeXchi\chi                         % new \chi, exactly on the baseline
\newbox\chibox
\setbox0 \hbox{\mathsurround0pt $\TeXchi$}
\setbox\chibox \hbox{\raise\dp0 \box 0 }
\def\chi{\copy\chibox}

% quadratino di fine dimostrazione

% abbreviazioni specifiche del lavoro

\def\Vp{V^*}

\def\normaV #1{\norma{#1}_V}

\def\kaz{\ka_*}
\def\kau{\ka^*}
\def\az{a_0}
\def\bz{b_0}

\def\hn{\hat n_\tau}
\def\hs{\hat s_\tau}
\def\hi{\hat i_\tau}
\def\hh{\hat h_\tau}

\def\on{\overline n_\tau}
\def\os{\overline s_\tau}
\def\oi{\overline i_\tau}
\def\oh{\overline h_\tau}

\def\un{\underline n_\tau}
\def\us{\underline s_\tau}
\def\ui{\underline i_\tau}
\def\uh{\underline h_\tau}

\def\nk{n_k}
\def\sk{s_k}
\def\ik{i_k}
\def\hk{h_k}
\def\ek{e_k}
\def\lak{\lambda_k}
\def\lakp{\lambda_{k+1}}

\def\nkp{n_{k+1}}
\def\skp{s_{k+1}}
\def\ikp{i_{k+1}}
\def\hkp{h_{k+1}}
\def\ekp{e_{k+1}}

\def\nm{n_m}

\def\soluz{(n,s,i,h)}

\def\nz{n^0}
\def\sz{s^0}
\def\iz{i^0}
\def\hz{h^0}
\def\nzt{\nz_\tau}
\def\szt{\sz_\tau}
\def\izt{\iz_\tau}
\def\hzt{\hz_\tau}
\def\utau{u_\tau}

\def\nmax{n^*}
\def\nmin{n_*}
\def\smax{s^*}
\def\imax{i^*}
\def\hmax{h^*}
\def\Tr{\tilde T}
\def\kat{\tilde\ka}
\def\At{\tilde A}
\def\tauz{\tau_0}

\def\signeps{\sign_\eps}
\def\modeps #1{|#1|_\eps}

\def\gn{g_{n,k}}
\def\gs{g_{s,k}}
\def\gi{g_{i,k}}
\def\gh{g_{h,k}}

%%%%%%%%%%%%%%%%%%%%%%%%%%%%%%
\Begin{document}
%%%%%%%%%%%%%%%%%%%%%%%%%%%%%%%%%

%%%%%%%%%%%%%%%%%%%%%%%%%%%%%%%%%
%% front page
%%%%%%%%%%%%%%%%%%%%%%%%%%%%%%%%%
%
\title{Well-posedness for a diffusion-reaction\\ compartmental model\\ simulating the spread of COVID-19}
\author{}
\date{}
\maketitle
\Bcenter
\vskip-1cm
{\large\sc Ferdinando Auricchio$^{(\circ)}$}\\
{\normalsize e-mail: {\tt ferdinando.auricchio@unipv.it}}\\[.25cm]
{\large\sc Pierluigi Colli$^{(*)}$}\\
{\normalsize e-mail: {\tt pierluigi.colli@unipv.it}}\\[.25cm]
{\large\sc Gianni Gilardi$^{(*)}$}\\
{\normalsize e-mail: {\tt gianni.gilardi@unipv.it}}\\[.25cm]
{\large\sc Alessandro Reali$^{(\circ)}$}\\
{\normalsize e-mail: {\tt alessandro.reali@unipv.it}}\\[.25cm]
{\large\sc Elisabetta Rocca$^{(*)}$}\\
{\normalsize e-mail: {\tt elisabetta.rocca@unipv.it}}\\[.45cm]
$^{(\circ)}$
{\small Dipartimento di Ingegneria Civile e Architettura, Universit\`a di Pavia}\\
{\small via Ferrata 3, 27100 Pavia, Italy}\\
{\small and Research Associate at the IMATI -- C.N.R. Pavia}\\[.25cm]
$^{(*)}$
{\small Dipartimento di Matematica ``F. Casorati'', Universit\`a di Pavia}\\
{\small via Ferrata 5, 27100 Pavia, Italy}\\
{\small and Research Associate at the IMATI -- C.N.R. Pavia}
\Ecenter
\Begin{abstract}
This paper is concerned with the well-posedness of a diffusion-reaction system for
a Susceptible-Exposed-Infected-Recovered (SEIR) mathematical model. 
This model
is written in terms of four nonlinear partial differential equations with nonlinear
diffusions, depending on the total amount of the SEIR populations. 
The model
aims at describing the spatio-temporal spread of the COVID-19 pandemic and is a
variation of the one recently introduced, discussed and tested in 
[A. Viguerie et al,
\pier{Diffusion-reaction} compartmental models \pier{formulated} in a continuum mechanics
framework: application to COVID-19\pier{,} mathematical analysis, and numerical study,
Comput.\ Mech.\ {\bf 66} (2020) 1131-1152]. 
Here, we deal with the mathematical analysis
of the resulting Cauchy--Neumann problem: 
the existence of solutions is proved in
a rather general setting and a suitable time discretization procedure is employed.
It is worth mentioning that the uniform boundedness of the discrete solution is shown
by carefully exploiting the structure of the system. 
Uniform estimates and passage
to the limit with respect to the time step allow to complete the existence proof.
Then, two uniqueness theorems are offered, one in the case of a constant diffusion
coefficient and the other for more regular data, in combination with a regularity
result for the solutions.
\vskip3mm


\End{abstract}

%%%%%%%%%%%%%%%%%%%%%%%%%%%%%%%%%%%%%%%%%%%%%%%%

\pagestyle{myheadings}
\newcommand\testopari{\sc Auricchio \ --- \ Colli \ --- \ Gilardi \ --- \ Reali \ --- \ Rocca}
\newcommand\testodispari{\sc Well-posedness for a model simulating the spread of COVID-19}
\markboth{\testodispari}{\testopari}

%%%%%%%%%%%%%%%%%%%%%%%%%%%%%%%%%
%% very beginning
%%%%%%%%%%%%%%%%%%%%%%%%%%%%%%%%%

\section{Introduction}
\label{Intro}
\setcounter{equation}{0}

\ale{The COVID-19 outbreak that deeply affected the world starting from the first months of 2020 led to a new, strong interest by researchers towards the development of mathematical models of infectious diseases, allowing the assessment of different scenarios and aiming at assisting the complex political decision-making process during the pandemic.
As a consequence, many papers were recently published, proposing interesting modeling ideas (see, e.g., \cite{Albi,Calleri,Gatto,Giordano,Jha,Linka,Parolini,Wang,Zohdi}), often based on compartmental models, where the considered population is divided into ``compartments'' based on their qualitative characteristics (like, e.g., ``susceptible'', ``infected'', ``recovered''), with different assumptions about the nature and rate of transfer across compartments. Despite this kind of models do not readily offer the possibility of a multiscale vision (as, e.g., proposed in \cite{Bellomo}), that would be a preferred feature given the nature of the phenomena to be simulated, they have the advantage of allowing a relatively easy introduction of diffusion terms (in this context, compelling ideas on diffusion models can be found, among others, in \cite{Bellomo2,Bellomo3,Winkler} and references therein). 
For a recent overview of mathematical models for virus pandemic, interested readers are referred to the articles included in \cite{Bellomo4}. Here, a compartmental model able to describe the spatial heterogeneity underlying the spread of an epidemic in a realistic city network is proposed in \cite{Bertaglia}. Another recent example where a classical compartmental model is enhanced to take into account infected individuals traveling on lines of fast diffusion is instead given by \cite{Berestycki}.
In fact, compartmental models are typically governed by a system of ordinary differential equations (ODEs) in time, which might possibly be enriched with coupling terms, additional equations, or other approaches to take into account also the spatial variation of the studied problem. An interesting alternative is to resort to compartmental models based on partial differential equations (PDEs), to accurately account for spatial variations at the continuum level.}

\ale{Following this latter idea, in the present contribution, we consider the PDE-based model introduced}
\oldgianni{in \cite{Vig1, Vig2} \pcol{(a variation of it with delay terms is discussed in \cite{Guglielmi})}.
In particular, the following system has been considered} 
\Bsist
  && \dt s
  = \alpha n - (1 - A_0/n) \beta_i si - (1 - A_0/n) \beta_e se - \mu s + \div(n \nu_s \nabla s) 
  \label{Ieqs}
  \\
  && \dt e
  = (1 - A_0/n) \beta_i si + (1 - A_0/n) \beta_e se - \sigma e - \phi_e e - \mu e + \div(n \nu_e \nabla e) 
  \qquad
  \label{Ieqe}
  \\
  && \dt i
  = \sigma e - \phi_d i - \phi_r i - \mu i + \div(n \nu_i \nabla i) 
  \label{Ieqi}
  \\
  && \dt r = \phi_r i + \phi_e e - \mu r + \div(n \nu_r \nabla r) 
  \label{Ieqr}
  \\
  && \dt d = \phi_d i \,.
  \label{Ieqd}
\Esist
\oldgianni{In the above equations, the symbols $s$, $e$, $i$, $r$, and $d$ denote
the susceptible population, the exposed population, the infected population, the recovered population, 
and the deceased population, respectively,
and $n:=s+e+i+r$ is the living population.
Note that $d$ refers only to deaths due to COVID-19. 
Due to the names of the compartments used, this model may be called a susceptible-exposed-infected-recovered-deceased (SEIRD) model.
Moreover, equations \accorpa{Ieqs}{Ieqr} are complemented by initial conditions
and no-flux boundary conditions,
while just an initial condition is associated to~\eqref{Ieqd}.}

\oldgianni{In \cite{Vig1} the authors presented a SEIRD mathematical model based on partial differential equations
coupled with a heterogeneous diffusion model. 
The model was used to describe the spatio-temporal spread of the COVID-19 pandemic 
and to capture dynamics also based on human habits and geographical features. 
To test the model, the outputs generated by a finite-element solver 
were compared with measured data over the Italian region of Lombardy 
and the results showed a strong qualitative agreement 
between the simulated forecast of the spatio-temporal COVID-19 spread in Lombardy 
and epidemiological data collected at the municipality level.} 

\oldgianni{In \cite{Vig2} the authors proposed a formulation of compartmental models 
based on partial differential equations through \pier{familiar} continuum mechanics concepts, 
interpreting such models in terms of fundamental equations of balance and compatibility, 
joined by a constitutive relation. 
Such an interpretation was useful to aid understanding and interdisciplinary collaboration. 
The authors formally derived the model sensitivity to diffusion, 
described its growth and decay, and established its stability in the $L^1$ norm. 
Attention was paid to an ODE version of the model, 
using it to derive a basic reproduction number $R_0$ as well as analyzing its spectrum. 
Additionally, a series of numerical simulations showed the role that numerical methods, diffusion, 
and specific ingredients played in the behavior of the system. 
Implicit models were found to be effective in describing the temporal dynamics of the system, 
and \pier{second-order} in-time methods in particular.}
\ale{The model proposed in \cite{Vig2} and used for simulating the COVID-19 spread in Lombardy, was further exploited in \cite{Grave,Grave2}. In particular, finite element simulations were carried out in the context of adaptive mesh refinement and coarsening, forecasting the COVID-19 spread also in the U.S. state of Georgia and in the Brazilian state of Rio de Janeiro. Good agreements with real-world epidemiological data were attained in both time and space.}

\oldgianni{Coming back to system \accorpa{Ieqs}{Ieqd},}
we note that $d$ can be determined from \eqref{Ieqd} and the associated initial condition
whenever the system of the first four equation in solved.
Thus, since we are just interested in well-posedness,
we do not consider the last equation
\oldgianni{and study system \accorpa{Ieqs}{Ieqr}.
Namely, we consider a slightly modified system.
Indeed, we replace the terms $1-A_0/n$ and $\phi_d i$ appearing in \accorpa{Ieqs}{Ieqe} and \eqref{Ieqi}
by a non-singular one of the form~$A(n)$ and by the product~$\phi_d in$, respectively.
On the contrary, we can treat the degeneracy in the diffusion terms
that occur as $n$ approaches zero, by proving that $n$ is bounded away from zero. 
This is done even in the more general situation of a coefficient of the form~$\ka(n)$,
where $\ka$ is continuous and strictly positive on~$(0,+\infty)$.
Finally, we assume that the constants $\nu_s$, $\nu_e$, $\nu_i$ and $\nu_r$ are the same.
For this modified model, we prove an existence result under mild assumption on the initial data.
Furthermore, we give two uniqueness results:
the main assumption for the first one is that $\ka$ is a positive constant,
while the second one requires that the nonlinearities and the initial data are smooth.
Under the latter assumptions, we prove a regularity result as well.} 

\pier{We hope that our results could capture the interest for the resulting systems and give rise to new approximations and simulations, in order to deal with and test our variations in the framework of some geographical areas. Moreover, we claim that our arguments for well-posedness provide a serious mathematical validation to this class of SEIR models.}

The paper is organized as follows. 
In the next section, we list our assumptions and notations
and state our results.
The proof of \oldgianni{Theorem~\ref{Existence}} regarding the existence of a solution 
is given in Section~\ref{EXISTENCE}
and is prepared in Section~\ref{DISCRETE} where an approximating problem
obtained by time discretization is studied.
Finally, Section~\ref{UNIQUENESS} is devoted to the uniqueness of the solution.

%%%%%%%%%%%%%%%%%%%%%%%%%%%%%%%%%%%%%%%%%%%%%%%%%%%%%%%%%%%%%%%%%%%%%%%%

\section{Statement of the problem and results}
\label{STATEMENT}
\setcounter{equation}{0}

In this section, we state precise assumptions and notations and present our results.
First of all, the subset $\Omega\subset\erre^3$
\oldgianni{(lower-dimensional cases can be treated in the same way)} 
is~assumed to be bounded, connected and smooth.
The symbol $\dn$ denotes for the normal derivative on $\Gamma:=\partial\Omega$.
Moreover, we set for brevity
\Beq
  Q_t := \Omega \times (0,t)
  \aand
  Q := \Omega \times (0,T).
  \label{defQt}
\Eeq
If $X$ is a Banach space, $\norma\cpto_X$ denotes both its norm and the norm of~$X^d$.
The only exception from this convention on the norms is given
by the spaces $L^p$ ($1\leq p\leq\infty$) constructed on $(0,T)$,
$\Omega$ and~$Q$, 
whose norms are often denoted by~$\norma\cpto_p$,
and by the space $H$ defined below and its powers,  
whose norms are simply denoted by~$\norma\cpto$.
We~put
\Bsist
  && H := \Ldue \,, \quad  
  V := \Huno 
  \aand
  W := \graffe{v \in \Hdue: \ \dn v = 0} .
  \label{defspazi}
\Esist 
Moreover, $\Vp$~is the dual space of $V$ and $\<\cpto,\cpto>$ is the dual pairing between $\Vp$ and~$V$.
In the sequel, we work in the framework of the Hilbert triplet
$(V,H,\Vp)$ \oldgianni{obtained by identifying $H$ with a subspace of $\Vp$ in the usual way}.
Thus, by also using the symbol $(\cpto,\cpto)$ for the standard inner product of~$H$
(this symbol will be used for any power of $H$ as well later~on),
we have
$\<g,v>=(g,v)$
for every $g\in H$ and $v\in V$.

Now, we list our assumptions on the structure of our system. 
We just consider the functions $A$ and $\ka$ mentioned in the introduction
and the constants $\alpha$ and $\mu$,
since all the other constants appearing in \accorpa{Ieqs}{Ieqr} will be normalized in the following.
However, we keep them for a while, for the reader's convenience.
We assume~that
\Bsist
  && \hbox{$A,\ka:(0,+\infty)\to\erre$ are continuous, with }
  \non
  \\
  && \quad \hbox{$A$ nonnegative and $\ka$ strictly positive}
  \label{hpAk}
  \\
  && \hbox{$\alpha$ and $\mu$ are positive constants}
  \label{hpam}
\Esist
\Accorpa\HPstruttura hpAk hpam
and we observe that a term like $\ka(n)=n$ is allowed in the equations.
In principle, the system we are interested in is the following 
\Bsist
  && \dt s + A(n) \beta_i si + A(n) \beta_e se + \mu s - \div(\ka(n)\nabla s)
  = \alpha n 
  \label{eqs}
  \\
  && \dt e -  A(n) \beta_i si -  A(n) \beta_e se + \sigma e + \phi_e e + \mu e - \div(\ka(n)\nabla e) 
  = 0
  \qquad
  \label{eqe}
  \\
  && \dt i + \oldgianni{\phi_d in} + \phi_r i + \mu i - \div(\ka(n)\nabla i) 
  = \sigma e
  \label{eqi}
  \\
  && \dt r - \phi_r i - \phi_e e + \mu r - \div(\ka(n)\nabla\oldgianni r) 
  = 0 
  \label{eqr}
\Esist
where each equation is complemented by no-flux boundary conditions and an initial condition.
However, it is more convenient to consider an equivalent system in the unknowns
$n$, $s$, $i$ and~$h$, where $n$ and $h$ are related to the \oldgianni{other} unknowns~by
\Beq
  n := s + e + i + r
  \aand
  h := s + e \,.
  \label{defnh}
\Eeq
The new system is obtained from the previous one by summing up all the equations,
keeping \eqref{eqs} and \eqref{eqi} and adding \eqref{eqe} to~\eqref{eqs}.
Hence, it is given~by 
\Bsist
  && \dt n + \phi_d i n - \div(\ka(n)\nabla n) 
  = (\alpha-\mu) n
  \label{eqprima}
  \\
  && \dt s + A(n) \beta_i si + A(n) \beta_e s (h-s) + \mu s - \div(\ka(n)\nabla s)
  = \alpha n 
  \label{eqseconda}
  \\
  && \dt i + \oldgianni{\phi_d in} + \phi_r i + \mu i - \div(\ka(n)\nabla i) 
  = \sigma (h-s)
  \label{eqterza}
  \\
  && \dt h + \mu h +(\sigma+\phi_e) h - \div(\ka(n)\nabla h) 
  = \alpha n + (\sigma+\phi_e) s \,.
  \label{eqquarta}
\Esist
However, as said above, among all of the constants appearing in these equations,
just $\alpha$ and $\mu$ will play a significant role in the mathematical treatment.
Hence, we replace the other constants and even some sums of them by~$1$, without loss of generality.
At this point, we are ready to give our notion of solution.
For the initial data we require that 
\Beq
  \nz \,, \sz \,, \iz \,, \hz \in \Linfty
  \quad \hbox{satisfy} \quad
  \inf\nz > 0 , \quad
  \hz \geq \sz \geq 0
  \aand
  \iz \geq 0
  \label{hpdati}
\Eeq
and we call solution a quadruplet $\soluz$ enjoying the requirements
\Bsist
  && n \,, s \,, i \,, h \in \H1\Vp \cap \C0 H \cap \L2V \cap \LQ\infty
  \qquad
  \label{regsoluz}
  \\
  && \inf n > 0 ,\quad
  h \geq s \geq 0 
  \aand
  i \geq 0
  \quad \aeQ
  \label{segni}
\Esist
\Accorpa\Regsoluz regsoluz segni
and satisfying, \aet\ and for every $v\in V$, the variational equations 
\Bsist
  && \< \dt n , v >
  + \iO i n v
  + \iO \ka(n) \nabla n \cdot \nabla v
  = (\alpha-\mu) \iO n v
  \label{prima}
  \\
  && \< \dt s , v >
  + \iO \bigl( A(n)si + A(n) s (h-s) + \mu s \bigr) v
  + \iO \ka(n) \nabla s \cdot \nabla v
  \non
  \\
  && \quad {} = \alpha \iO n v
  \label{seconda}
  \\
  \separa
  && \< \dt i , v >
  + (1+\mu) \iO i v
  + \iO \ka(n) \nabla i \cdot \nabla v
  = \iO (h-s) v
  \label{terza}
  \\
  && \< \dt h , v >
  + \oldgianni{\iO n i v}
  + (1+\mu) \iO h v
  + \iO \ka(n) \nabla h \cdot \nabla v
  = \iO (\alpha n + s) v
  \qquad
  \label{quarta}
\Esist
and the initial conditions
\Beq
  (n,s,i,h)(0) = (\nz,\sz,\iz,\hz) \,.
  \label{cauchy}
\Eeq
\Accorpa\Pbl prima cauchy

\Brem
\label{IntPbl}
We notice that the regularity $\C0H$ explicitly stated in \eqref{regsoluz}
in fact follows from the other conditions since
$\H1\Vp\cap\L2V\subset\C0H$.
So we do not mind \oldgianni{about} it in the following.
We also observe that the above variational equations are equivalent 
to their integrated versions with time dependent test functions.
For instance, \eqref{prima} is equivalent~to
\Bsist
  && \ioT \< \dt n(t) , v(t) > \, dt
  + \intQ i n v
  + \intQ \ka(n) \nabla n \cdot \nabla v
  = (\alpha-\mu) \intQ n v
  \non
  \\
  && \quad \hbox{for every $v\in\L2V$}.
  \label{intprima}
\Esist
\Erem

Here is our \oldgianni{first} result.
\Bthm
\label{Existence}
Assume \HPstruttura\ and \eqref{hpdati}.
Then, there exists at least one quadruplet $\soluz$
that satisfies the conditions \Regsoluz\ and solves problem \Pbl.
\Ethm

As for uniqueness, we present two results.
The first one regards a particular case and its proof, given in Section~\ref{UNIQUENESS1}, is elementary.

\Bthm
\label{Uniqueness1}
In addition to \HPstruttura\ and \eqref{hpdati},
assume that $A$ is locally \Lip\ continuous and that $\ka$ is a \oldgianni{positive} constant.
Then, the solution to problem \Pbl\ satisfying \Regsoluz\ is unique.
\Ethm

Our last result regards uniqueness in the case of a non-constant~$\ka$.
Its proof, given in Section~\ref{UNIQUENESS2}, is much more involved and  
is strictly related to a high regularity of the solution.
For this reason, we have to assume that both the nonlinearities and the initial data are much smoother.

\Bthm
\label{Uniqueness2}
In addition to \HPstruttura\ and \eqref{hpdati},
assume that
\Bsist
  && \hbox{$A$ is locally \Lip\ continuous and $\ka$ is a $C^1$ function}
  \label{hpnonlin}
  \\
  && \nz,\, \sz,\, \iz,\, \hz \in \Wx{2,\infty} 
  \quad \hbox{with zero normal derivatives on $\Gamma$}.
  \label{hpdatareg}
\Esist
Then, the solution to problem \Pbl\ satisfying \Regsoluz\ is unique
and enjoys the following regularity properties
\Beq
  n,\, s,\, i,\, h \in \W{1,p}{\Lx p} \cap \L p{\Wx{2,p}}
  \quad \hbox{for every $p\in[1,+\infty)$}.
  \label{highreg}
\Eeq
\Ethm

Throughout the paper, we make use of
the \Holder\ inequality and the Sobolev inequality 
related to the continuous embedding $V\subset L^p(\Omega)$ with $p\in[1,6]$
(since $\Omega$ is three-dimensional bounded and smooth).
We also account for the elementary identity and inequalities
\Bsist
  \hskip-1cm&& a (a-b)
  = \frac 12 \, a^2
  + \frac 12 \, (a-b)^2
  - \frac 12 \, b^2
  \geq \frac 12 \, a^2
  - \frac 12 \, b^2
  \quad \hbox{for every $a,b\in\erre$},
  \label{elementare}
  \\
  \hskip-1cm&& ab\leq \delta a^2 + \frac 1 {4\delta}\,b^2
  \quad \hbox{for every $a,b\in\erre$ and $\delta>0$},
  \label{young}
\Esist
\Accorpa\Elementari elementare young
and quote \eqref{young} as the Young inequality.
Furthermore, we take advantage of the summation by parts formula
\Beq
  \somma k0{m-1} a_{k+1} (b_{k+1} - b_k)
  = a_m b_m - a_1 b_0
  - \somma k1{m-1} (a_{k+1} - a_k) b_k\,,
  \label{byparts}
\Eeq
which is valid for arbitrary real numbers $a_1,\dots,a_m$ and $b_0,\dots,b_m$.

%%%%%%%%%%%%%%%%%%%%%%%%%%%%%%%%%%%%%%%%%%%%%%%%%%%%%%%%%%%%%%%%%%%%%%%%

\section{The discrete problem}
\label{DISCRETE}
\setcounter{equation}{0}

In this section, we prepare \oldgianni{the proof of Theorem~\ref{Existence}}
by introducing and solving an approximating problem obtained by time discretization.
However, the structural functions $A$ and $\ka$ have to satisfy different assumptions
 and the initial data have to be smoother.
In the next section, by starting from the original structure and the original initial data,
we consider the discrete problem with structural functions and approximating initial data
constructed in order to satisfy the assumptions listed below.
\oldgianni{Two~constants $\kaz$ and $\kau$} and two real functions $\At$ and $\kat$ defined in the whole of $\erre$ 
are given such~that
\Bsist
  && \oldgianni{\kau \geq \kaz > 0}
  \aand
  \hbox{$\At,\kat:\erre\to\erre$ are continuous with}
  \non
  \\
  && \quad
  \At(y) \geq 0 
  \aand
  \oldgianni{\kaz \leq \kat(y) \leq \kau}
  \quad \hbox{for every $y\in\erre$}.
  \label{hpAkdiscr}
\Esist

\Bnot
\label{Akatilde}
However, we prefer to use the lighter symbols $A$ and $\ka$ 
instead of the heavy $\At$ and~$\kat$.
Indeed, no confusion can arise since the original functions 
$A$ and $\ka$ introduced in \eqref{hpAk} will never appear within the section.
\Enot
For a fixed positive integer~$N$, we set $\tau:=T/N$.
Then, the time-discretized problem we are going to study is the following:
given 
\Beq
  \nzt \,, \szt \,, \izt \,, \hzt \in V\cap\Linfty
  \quad \hbox{with} \quad
  \nzt \geq 0 , \quad
  \hzt \geq \szt \geq 0
  \aand
  \izt \geq 0
  \label{hpdatiV}
\Eeq
we look for four $(N+1)$-tuples 
\Bsist
  && (n_0,n_1,\dots,n_N), \
  (s_0,s_1,\dots,s_N), \
  (i_0,i_1,\dots,i_N), \
  (h_0,h_1,\dots,h_N) 
  \non
  \\
  && \quad \hbox{belonging to $(V\cap\Linfty)^{N+1}$}
  \label{tuples}
\Esist
satisfying, for $k=0,\dots,N-1$, the variational equations
\Bsist
  && \bigl( \frac {\nkp-\nk}\tau , v \bigr)
  + (\ik \nkp , v)
  + (\ka(\nkp) \nabla\nkp , \nabla v)
  = (\alpha-\mu) (\nkp , v)
  \label{primak}
  \\
  && \bigl( \frac {\skp-\sk}\tau , v \bigr)
  + \bigl( A(\nkp)\skp\ik + A(\nkp)\skp(\hk-\sk) + \mu\skp , v \bigr)
  \non
  \\
  && \quad {}
  + \bigl( \ka(\nkp) \nabla\skp , \nabla v \bigr)
  = \alpha (\nkp , v)
  \label{secondak}
  \\
  && \bigl( \frac {\ikp-\ik}\tau , v \bigr)
  \oldgianni{{}+(\nkp\ikp , v)}
  + (1+\mu) (\ikp , v)
  + \bigl( \ka(\nkp)\nabla\ikp , \nabla v \bigr)
  \non
  \\
  && \quad {}
  = ( \hkp-\skp , v)
  \qquad
  \label{terzak}
  \\
  && \bigl( \frac {\hkp-\hk}\tau , v \bigr)
  + (1+\mu) (\hkp , v)
  + \bigl( \ka(\nkp)\nabla\hkp , \nabla v \bigr)
  \non
  \\
  && \quad {}
  = (\alpha\nkp + \skp , v)
  \label{quartak}
\Esist
all for every $v\in V$,
as well as the sign and initial conditions
\Bsist
  && \nk \geq 0 ,\quad
  \hk \geq \sk \geq 0 
  \aand
  \ik \geq 0
  \quad \hbox{for $k=1,\dots,N$}
  \label{segnik}
  \\
  && n_0 = \nzt \,, \quad
  s_0 = \szt \,, \quad
  i_0 = \izt 
  \aand
  h_0 = \hzt .
  \label{discrcauchy}
\Esist
\Accorpa\Pbltau primak discrcauchy

\Bthm
\label{Wellposednessdiscr}
Under assumptions \eqref{hpAkdiscr} and \eqref{hpam} on the structure,
suppose that
\Beq
  \frac 1\tau - \alpha + \mu > 0 .
  \label{hptau}
\Eeq
Then, if the initial data satisfy \eqref{hpdatiV}, problem \Pbltau\ has a unique solution.
\Ethm

We prepare an easy lemma.

\Blem
\label{Gianni}
Let $\az$ and $\bz$ be positive constants and let $a,b,f\in\Linfty$ satisfy 
\Beq
  a \geq \az \,, \quad b \geq \bz \aand f \geq 0 \quad \aeO .
  \label{hpG}
\Eeq
Then, the problem of finding $u\in V$ satisfying the variational problem
\Beq
  \iO a \nabla u \cdot \nabla v
  + \iO b u v 
  = \iO f v
  \quad \hbox{for every $v\in V$}
  \label{pblG}
\Eeq
has a unique solution, and this solution satisfies
\Beq
  0 \leq u \leq \bz^{-1} \norma f_\infty \quad \aeO.
  \label{tesiG}
\Eeq
Moreover, if $\lambda\in\erre$ and $f\geq\lambda\,b$ \aeO, then $u\geq\lambda$ \aeO.
\Elem

\Bdim
The existence of a unique solution is clear 
since the \lhs\ of \eqref{pblG} is an inner product in $V$ that is equivalent to the usual one.
Moreover, the first inequality of \eqref{tesiG} is given by the weak maximum principle
and the same can be said for the last sentence of the statement, 
since the equation obtained by replacing $f $ by $f-\lambda\,b$ is satisfied by~$u-\lambda\,$.
Now, we prove the second inequality of~\eqref{tesiG}.
\oldgianni{%
We set $w:=u-\bz^{-1}\norma f_\infty$. 
Then, we have for every $v\in V$
\Bsist
  && \iO a \nabla w \cdot \nabla v
  + \iO b w v
  = \iO a \nabla u \cdot \nabla v
  + \iO b u v
  - \iO b \, \bz^{-1} \norma f_\infty v
  \non
  \\
  && = \iO \bigl( f - b \, \bz^{-1} \norma f_\infty \bigr) v.
  \non
\Esist
Since $f-b\,\bz^{-1}\norma f_\infty\leq f-\norma f_\infty\leq0$,
we have $w\leq0$ by the weak maximum principle, 
i.e., the desired inequality.}
\Edim

\step
Proof of Theorem~\ref{Wellposednessdiscr}

It suffices to prove the following:
for $k=0,\dots,N-1$, if $(\nk,\sk,\ik,\hk)$ belongs to $(V\cap\Linfty)^4$ 
and satisfies the inequalities in \eqref{segnik},
then system \accorpa{primak}{quartak} has a unique solution $(\nkp,\skp,\ikp,\hkp)$
belonging to $(V\cap\Linfty)^4$ and satisfying the analogous inequalities,
i.e., $\nkp\geq0$, $\hkp\geq\skp\geq0$ and $\ikp\geq0$.
We recall \eqref{hpdatiV} for the case $k=0$,
fix $k$ and $(\nk,\sk,\ik,\hk)$ as said 
and show that we can find a unique solution to \accorpa{primak}{quartak} with the proper sign conditions.
This is done in two steps.

\step
Solution to the first equation

We introduce the function $K:\erre\to\erre$ by setting
\Beq
  K(y) := \int_0^y \ka(z) \, dz
  \quad \hbox{for $y\in\erre$}
  \label{defK}
\Eeq
and we would like to assume $u:=K(\nkp)$ as the new unknown for~\eqref{primak}.
Due to \eqref{hpAkdiscr} for~$\ka$,
the function $K$ is one-to-one, onto and \Lip\ continuous, and its inverse is \Lip\ continuous too.
Thus, equation \eqref{primak} can be rewritten in term of~$u$.
\gianni{Namely, noting that $\nabla u=\nabla(K(\nkp))=\ka(\nkp)\nabla\nkp$,
we can write it as the variational formulation of the homogeneous Neumann problem for the equation
\Beq
  \lambda K^{-1}(u)
  - \Delta u 
  = (1/\tau)\nk
  \quad \hbox{where} \quad
  \lambda := \frac 1\tau + \ik - \alpha + \mu \,.
  \label{barbu}
\Eeq
Notice that every variational solution automatically belongs to $W$ 
and satisfies both \eqref{barbu} and the homogeneous Neumann condition due to elliptic regularity.
We claim that the new problem has a unique solution.
Indeed, it is equivalent to the minimization of the functional $J:V\to\erre$ given~by
\Beq
  J(v)
  := \frac 12 \iO |\nabla v|^2
  + \iO \lambda \, \calK(v)
  - \frac 1\tau \iO \nk \, v 
  \label{minimo}
\Eeq
where 
\Beq
  \calK(r) := \int_0^r K^{-1}(s) \, ds
  \quad \hbox{for $r\in\erre$}
  \non
\Eeq
and we show that this problem has a unique solution.
Notice that $J$ actually is well-defined, 
since $K^{-1}$ is \Lip\ continuous, $\lambda\in\Linfty$ and $\nk\in H$. 
Moreover, $J$~is strictly convex and coercive.
To see this, it suffices to recall that $\inf\lambda>0$ 
(since we are assuming \eqref{hptau} and $\ik$ is nonnegative)
and observe that
$\calK''=(K^{-1})'=1/(\ka\circ K^{-1})\geq1/\kau$.
Hence, the above problem has a unique solution $u$
and we conclude that $\nkp=K(u)$ is the unique solution to~\eqref{primak}.
Since $u\in W$, we deduce that $u$ is bounded
so that that $\nkp$ is bounded too.
To see that it is nonnegative, it suffices to rearrange \eqref{primak} and apply Lemma~\ref{Gianni}.}

\step
Solution to the other equations and conclusion

Since $\nkp$ is known, we can solve \eqref{secondak} for~$\skp$
by applying the first part of Lemma~\ref{Gianni}.
We can do the same to find $\hkp$ and $\ikp$ from \eqref{quartak} and \eqref{terzak} in this order.
The lemma also ensures that $\skp$, $\hkp$ and $\ikp$ are bounded and nonnegative,
\oldgianni{provided we can prove that} $\hkp\geq\skp$.
To this end, we set $\ek:=\hk-\sk$ and $\ekp=\hkp-\skp$, 
take the difference between \oldgianni{\eqref{quartak} and \eqref{secondak}}
and write it in the form
\Bsist
  && \iO \ka(\nkp) \nabla\ekp \cdot \nabla v
  + ((1/\tau)+1+\mu) \iO \ekp v
  \non
  \\
  && = \iO \bigl( A(\nkp)\skp\ik + A(\nkp)\skp\ek + (1/\tau)\ek \bigr) v
  \non
\Esist
for every $v\in V$.
Then, Lemma~\ref{Gianni} yields that $\ekp\geq0$.
This completes the proof.

%%%%%%%%%%%%%%%%%%%%%%%%%%%%%%%%%%%%%%%%%%%%%%%%%%%%%%%%%%%%%%%%%%%%%%%%

\section{Existence}
\label{EXISTENCE}
\setcounter{equation}{0}

This section is devoted to the proof of Theorem~\ref{Existence}.
Our argument relies on a priori estimates on the solution to a suitably specified discrete problem
and on the convergence of proper interpolating functions.
So, we introduce some notations concerning interpolation at once.

\Bnot
\label{Interpol}
Let $N$ be a positive integer and $Z$ be a Banach space.
We set $\tau:=T/N$ and $I_k:=((k-1)\tau,k\tau)$ for $k=1,\dots,N$.
Given $z=(z_0,z_1,\dots ,z_N)\in Z^{N+1}$,
we define the piecewise constant and piecewise linear interpolating functions
\Beq
  \overline z_\tau \in \L\infty Z , \quad
  \underline z_\tau \in \L\infty Z 
  \aand
  \hat z_\tau \in \W{1,\infty}Z
  \non
\Eeq
by setting 
\Bsist
  && \hskip -2em
  \overline z_\tau(t) = z^k
  \aand
  \underline z_\tau(t) = z^{k-1}
  \quad \hbox{for a.a.\ $t\in I_k$, \ $k=1,\dots,N$},
  \label{pwconstant}
  \\
  && \hskip -2em
  \hat z_\tau(0) = z_0
  \aand
  \dt\hat z_\tau(t) = \frac {z^k-z^{k-1}} \tau
  \quad \hbox{for a.a.\ $t\in I_k$, \ $k=1,\dots,N$}.
  \qquad
  \label{pwlinear}
\Esist
\Enot

For the reader's convenience,
we summarize the relations between the finite set of values
and the interpolating functions in the following proposition,
whose proof follows from \sfw\ computation:

\Bprop
\label{Propinterp}
With Notation~\ref{Interpol}, we have that
\Bsist
  && \norma{\overline z_\tau}_{\L\infty Z}
  = \max_{k=1,\dots,N} \norma{z^k}_Z \,, \quad
  \norma{\underline z_\tau}_{\L\infty Z}
   = \max_{k=0,\dots,N-1} \norma{z^k}_Z\,,
  \label{ouLinftyZ}
  \\
  && \norma{\dt\hat z_\tau}_{\L\infty Z}
  = \max_{0\leq k\leq N-1} \norma{(z^{k+1}-z^k)/\tau}_Z\,,
  \label{dtzLinftyZ}
  \\
  \separa
  && \norma{\overline z_\tau}_{\L2Z}^2
  = \tau \somma k1N \norma{z^k}_Z^2 \,, \quad
  \norma{\underline z_\tau}_{\L2Z}^2
  = \tau \somma k0{N-1} \norma{z^k}_Z^2 \,,
  \label{ouLdueZ}
  \\
  \separa
  && \norma{\dt\hat z_\tau}_{\L2Z}^2
  = \tau \somma k0{N-1} \norma{(z^{k+1}-z^k)/\tau}_Z^2\,, 
  \label{dtzLdueZ}
  \\
  && \norma{\hat z_\tau}_{\L\infty Z}
  = \max_{k=1,\dots,N} \max\{\norma{z^{k-1}}_Z,\norma{z^k}_Z\}
  = \max\{\norma{z_0}_Z,\norma{\overline z_\tau}_{\L\infty Z}\}\,,
  \qquad\qquad
  \label{hzLinftyZ}
  \\
  && \norma{\hat z_\tau}_{\L2Z}^2
  \leq \tau \somma k1N \bigl( \norma{z^{k-1}}_Z^2 + \norma{z^k}_Z^2 \bigr)
  \leq \tau \norma{z_0}_Z^2
  + 2 \norma{\overline z_\tau}_{\L2Z}^2 \,.
  \label{hzLdueZ}
\Esist
Moreover, it holds that
\Bsist
  && \norma{\overline z_\tau-\hat z_\tau}_{\L\infty Z}
  = \max_{k=0,\dots,N-1} \norma{z^{k+1}-{z^k}}_Z
  = \tau \, \norma{\dt\hat z_\tau}_{\L\infty Z}\,,
  \qquad
  \label{diffLinfty}
  \\
  && \norma{\overline z_\tau-\hat z_\tau}_{\L2Z}^2
  = \frac \tau 3 \somma k0{N-1} \norma{z^{k+1}-z^k}_Z^2
  \non
  \\
  && \oldgianni{{}= \frac 13 \, \norma{\overline z_\tau - \underline z_\tau}_{\L2Z}^2}
  = \frac {\tau^2} 3 \, \norma{\dt\hat z_\tau}_{\L2Z}^2\,,
  \label{diffLdue}
\Esist
and similar identities for the difference $\underline z_\tau-\hat z_\tau$.
Finally, we have that
\Bsist
  && \tau \somma k0{N-1} \norma{(z^{k+1}-z^k)/\tau}_Z^2 
  \leq \norma{\dt z}_{\L2Z}^2
  \non
  \\
  && \quad \hbox{if $z\in\H1Z$\aand $z^k=z(k\tau)$ for $k=0,\dots,N$}.
  \label{interpH1Z}
\Esist
\Eprop

We are now ready to properly specify the discrete problem.
Namely, starting from $A$ and $\ka$ as in~\eqref{hpAk},
we introduce new functions $\At$ and $\kat$ obtained by a truncation operator $\Tr$ 
to be used in the discrete problem.
We~set
\Bsist
  && \nmax := e^{2T(\alpha-\mu)^+} \norma\nz_\infty \,, \quad
  \smax := \norma\sz_\infty + T\alpha\nmax
  \non
  \\
  && \quad
  \hmax := \norma\hz_\infty + T(\alpha\nmax+\smax) \quad
  \imax := \norma{i_0}_\infty + T(\hmax+\smax) \quad
  \non
  \\
  && \aand
  \nmin := e^{-T(\imax+(\mu-\alpha)^+)} \inf\nz
  \label{defnstar}
\Esist
and define $\At,\kat:\erre\to\erre$ by setting for $y\in\erre$
\Beq
  \At(y) = A(\Tr(y)) 
  \aand
  \kat(y) = \ka(\Tr(y))
  \quad \hbox{where} \quad
  \Tr(y) := \max\graffe{\nmin, \min\graffe{y,\nmax}} .
  \label{kaAstar}
\Eeq
Next, we approximate the initial data $\nz$, $\sz$, $\iz$ and $\hz$ as in \eqref{hpdati}
by smoother functions $\nzt$, $\szt$, $\izt$ and $\hzt$ satisfying
\Bsist
  && \nzt \in V \cap \Linfty 
  \aand
  0 \leq \nzt \leq \norma\nz_\infty
  \quad \aeO
  \label{hpdatitau}
  \\
  && \norma\nzt \leq \norma\nz 
  \aand
  \tau\normaV\nzt^2 \leq \norma\nz^2
  \label{stimadatitau}
  \\
  && \nzt \to \nz 
  \quad \hbox{strongly in $H$ as $\tau\searrow0$}
  \label{convdatitau}
  \\
\noalign{\noindent and the analogues for $\szt$, $\izt$ and~$\hzt$ as well as \medskip}
  && \hzt \geq \szt
  \aand
  \nzt \geq \inf\nz
  \quad \aeO.
  \label{nztmin}
\Esist
This can be done, e.g., by a singular perturbation argument.
Indeed, if $\tau\in(0,1)$ and $u\in\Linfty$ is nonnegative, 
the unique solution $\utau\in V$ to the variational problem
\Beq
  \iO \bigl( \utau v + \tau \nabla \utau \cdot \nabla v \bigr)
  = \iO u v
  \quad \hbox{for every $v\in V$}
  \non
\Eeq
\oldgianni{belongs to $W\subset\Linfty$,}
satisfies 
$\inf u\leq\utau\leq\norma u_\infty$, $\norma\utau\leq\norma u$ 
and \oldgianni{$(1/2)\norma\utau^2+\tau\norma{\nabla\utau}^2\leq(1/2)\norma u^2$,
whence also $\tau\normaV\utau^2\leq\norma u^2$,}
and converges to $u$ strongly in $H$ as $\tau\searrow0$.
Finally, we assume 
\Beq
  \tau \in (0,\tauz)
  \quad \hbox{where} \quad
  \tauz \in (0,1)
  \aand
  \tauz \leq \frac 1 {2(\alpha-\mu)}
  \quad \hbox{if $\alpha>\mu$}.
  \label{hptauz}
\Eeq
The functions $\At$ and $\kat$ satisfy \eqref{hpAkdiscr} with
\Beq
  \oldgianni{\kaz := \min \graffe{\ka(y) : \ \nmin\leq y\leq\nmax}
  \aand
  \kau := \max \graffe{\ka(y) : \ \nmin\leq y\leq\nmax}.}
  \label{sceltakazkau}
\Eeq
Moreover, the requirements \eqref{hpdatiV} are fulfilled by the approximating data.
Furthermore, \eqref{hptauz} implies~\eqref{hptau}. 
Hence, we can solve the discretized problem 
with the specified structure and the approximating initial data.
We have the following basic result, 
which removes the supplementary assumptions~\eqref{hpAkdiscr}:

\Bprop
\label{Basic}
Under the restriction \eqref{hptauz}, let 
\Beq
  (n_0,n_1,\dots,n_N), \quad
  (s_0,s_1,\dots,s_N), \quad
  (i_0,i_1,\dots,i_N)
  \aand
  (h_0,h_1,\dots,h_N) 
  \label{soluztau}
\Eeq 
be as in \eqref{tuples} and solve the problem \Pbltau\ 
where one reads $\At$ and $\kat$ instead of $A$ and~$\ka$, respectively.
Then, \eqref{soluztau} solve the discrete problem with the original $A$ and~$\ka$ given in \eqref{hpAk}.
\Eprop

\Bdim
All the components $\nk$ are nonnegative.
However, as $A$ and $\ka$ are defined in the open interval $(0,+\infty)$, 
we have to reinforce this and show that each $\nk$ is strictly positive
in order to give a meaning to $A(\nk)$ and~$\ka(\nk)$.
We prove~that
\Beq
  \nmin \leq \nk \leq \nmax 
  \quad \hbox{\aeO \quad for $k=0,\dots,N$}.
  \label{basic}
\Eeq
Since $\At(y)=A(y)$ and $\kat(y)=\ka(y)$ for every $y\in[\nmin,\nmax]$,
this also yields that $\At(\nk)=A(\nk)$ and $\kat(\nk)=\ka(\nk)$ for $k=0,\dots,N$, 
i.e., the thesis of the statement.

The proof of \eqref{basic} is done in several steps.
The last of them needs $L^\infty$ bounds for the other components of the solution,
and these are proved as~well.

\step
First upper bound

We assume $0\leq k<N$ and check that we can apply Lemma~\ref{Gianni} to~$\nkp$.
Indeed, \eqref{primak} with $\kat$ in place of $\ka$ 
implies that $\nkp$ solves \eqref{pblG} with
\Beq
  a = \kat(\nkp) , \quad
  b = \frac 1\tau + \ik -\alpha + \mu
  \aand
  f = \frac 1\tau \, \nk \,.
  \label{gianninkp}
\Eeq
Moreover, we can take $\az=\kaz$ and $\bz=(1/\tau)-(\alpha-\mu)^+$ since $\ik$ is nonnegative.
Notice that $\bz>0$ by \eqref{hptauz}.
As $\nkp$ is nonnegative, the lemma implies that
\Beq
  \norma\nkp_\infty
  \leq \bz^{-1} \norma f_\infty
  = \frac \tau {1-\tau(\alpha-\mu)^+} \, \norma{(1/\tau)\nk}_\infty 
  = \frac 1 {1-\tau(\alpha-\mu)^+} \, \norma\nk_\infty \,.
  \non
\Eeq
As this holds for $0\leq k<N$, $n_0=\nzt$ and 
$\norma\nzt_\infty\leq\norma\nz_\infty$ by~\eqref{hpdatitau}, 
we deduce that
\Beq
  \norma\nk_\infty
  \leq (1-\tau(\alpha-\mu)^+)^{-k} \norma{\nzt}_\infty
  \leq (1-\tau(\alpha-\mu)^+)^{-N} \norma\nz_\infty
  \label{forbasic}
\Eeq
for $k=0,\dots,N$.
Now, observe that $\ln(1-y)\geq-2y$ for every $y\in(0,1/2]$.
Indeed, the logarithm is a concave function 
and the inequality is satisfied at $y=1/2$ and for $y>0$ small enough 
by comparison of the derivatives at~$0$.
Since $\tau(\alpha-\mu)^+\leq1/2$ by~\eqref{hptauz},
it follows that
\Beq
  \ln(1-\tau(\alpha-\mu)^+)
  \geq -2\tau(\alpha-\mu)^+ 
  \non
\Eeq
and \eqref{forbasic} implies
\Beq
  \norma\nk_\infty
  \leq (1-\tau(\alpha-\mu)^+)^{-N} \norma\nz_\infty
  \leq e^{2N\tau(\alpha-\mu)^+} \norma\nz_\infty
  = \nmax 
  \quad \oldgianni{\hbox{for $k=0,\dots,N$.}}
  \label{stimank}
\Eeq
Therefore, the second inequality of \eqref{basic} is proved.

Now, we prove analogous upper bounds for the other components of the solution
by applying Lemma~\ref{Gianni} once more, always with $a=\kat(\nkp)$.

\step
Further upper bounds

We start with $\sk$.
Equation \eqref{pblG} is solved by $\skp$ with
\Beq
  b = \frac 1\tau + \At(\nkp)\ik + \At(\nkp)(\hk-\sk) + \mu
  \aand 
  f = \frac 1\tau \sk + \alpha\nkp \,.
  \non
\Eeq
By recalling that $\At$ is nonnegative and that $\ik\geq0$ and $\hk\geq\sk$ \aeO, 
we can take $\bz=(1/\tau)+\mu$.
Since $\skp$ is nonnegative, the lemma yields
\Beq
  \norma\skp_\infty
  \leq \frac \tau {1+\tau\mu} \, \norma{(1/\tau)\sk+\alpha\nkp}_\infty
  \leq \norma\sk_\infty + \tau\alpha\norma\nkp_{\oldgianni\infty} \,.
  \non
\Eeq
On account of \oldgianni{\eqref{stimank} and \eqref{defnstar}}, we deduce that
\Beq
  \norma\sk_\infty
  \leq \norma{s_0}_\infty + k\tau\alpha\nmax
  \leq \smax 
  \quad \hbox{for $k=0,\dots,N$}.
  \label{sinfty}
\Eeq
For $\hkp$ we can take
\Beq
  b = \frac 1\tau + 1 + \mu
  \aand
  f = \frac 1\tau \, \hk + \alpha\,\nkp + \skp \,.
  \non
\Eeq
Then the lemma and the estimates for $\nk$ and $\sk$ already proved give
\Beq
  \norma\hkp_\infty 
  \leq \frac 1 {(1/\tau)+1+\mu} \, \norma{(1/\tau)\hk+\alpha\nkp+\skp}_\infty
  \leq \norma\hk_\infty + \tau (\alpha\nmax+\smax) \,.
  \non
\Eeq
Hence
\Beq
  \norma\hk_\infty
  \leq \norma{h_0}_\infty + k\tau (\alpha\nmax+\smax)
  \leq \hmax
  \quad \hbox{for $k=0,\dots,N$}.
  \label{hinfty}
\Eeq
Finally, for $\ikp$ we can take 
\Beq
  b = \frac 1\tau + \oldgianni\nkp + 1 + \mu
  \aand
  f = \frac 1\tau \, \ik + \hkp - \skp 
  \non
\Eeq
and a quite similar calculation also using \eqref{hinfty} yields
\Beq
  \norma\ik_\infty 
  \leq \imax
  \quad \hbox{for $k=0,\dots,N$}.
  \label{iinfty}
\Eeq

\step
Lower bound

Finally, we prove the \oldgianni{left} inequality of~\eqref{basic}.
We fix $k$ with $0\leq k<N$ and assume that $\nk\geq\lak$ \aeO\ for some $\lak>0$
and we prove that $\nkp\geq\lakp$ \aeO\ where
\Beq
  \lakp := \frac \lak {1+\tau(\imax+(\mu-\alpha)^+)} \,.
  \non
\Eeq
As already observed, \eqref{pblG} is solved by $\nkp$ with
$a$, $b$ and $f$ given by \eqref{gianninkp}.
Let us estimate $f-\lakp b$ from below.
We~have
\Beq
  f - \lakp b
  = \frac 1\tau \, \nk 
  - \frac \lak {1+\tau(\imax+(\mu-\alpha)^+)} \, \Bigl( \frac 1\tau + \ik + \oldgianni{\mu-\alpha} \Bigr)
  \geq \frac 1\tau \, \nk - \frac 1\tau \, \lak  
  \geq 0 \,.
  \non
\Eeq
Then, we can apply the last sentence of Lemma~\ref{Gianni} and obtain that
$\nkp\geq\lakp$ \aeO.
Since $n_0=\nzt\geq\inf\nz$ \aeO, we can take $\lambda_0=\inf\nz$ and~have
\Beq
  \nk \geq (\inf\nz) \, \bigl( 1+\tau(\imax+(\mu-\alpha)^+) \bigr)^{-k}
  \quad \aeO , \quad \hbox{for $k=0,\dots,N$}.
  \label{forlower}
\Eeq
Since $\ln(1+y)\leq y$ for every $y>-1$, we obtain
\Beq
  -k \ln \bigl( 1+\tau(\imax+(\mu-\alpha)^+) \bigr) 
  \geq -k \tau(\imax+(\mu-\alpha)^+)
  \geq -T(\imax+(\mu-\alpha)^+) 
  \non
\Eeq
whence also
\Beq
  \bigl( 1+\tau(\imax+(\mu-\alpha)^+) \bigr)^{-k}
  \geq e^{-T(\imax+(\mu-\alpha)^+)} .
  \non
\Eeq
Hence, \eqref{forlower} yields that $\nk\geq\nmin$ \aeO\ for every~$k$.
This concludes the proof.
\Edim

At this point, we establish two estimates whose proof is made very easy by the above $L^\infty$-bounds.
For the second one it is convenient to rewrite the estimates in terms of the interpolating functions.

\Bnot
\label{Constants}
For the sake of simplicity, in the rest of the section,
we use the same symbol $c$ for possibly different constants
that only depend on $\ka$, $A$, $\alpha$, $\mu$, $T$, $\Omega$ 
and the initial data appearing in~\eqref{hpdati}.
Hence, the values of $c$ do not depend on~$\tau$
and might change from line to line and even within the same line.
\Enot

\step
First a priori estimate

We rewrite \accorpa{primak}{quartak} in the form
\Bsist
  && \bigl( \frac {\nkp-\nk}\tau , v \bigr)
  + (\ka(\nkp) \nabla\nkp , \nabla v)
  + \kaz (\nkp , v)
  = (\gn , v)
  \label{primakbis}
  \\
  && \bigl( \frac {\skp-\sk}\tau , v \bigr)
  + \bigl( \ka(\nkp) \nabla\skp , \nabla v \bigr)
  + \kaz (\skp , v)
  = (\gs , v)
  \label{secondakbis}
  \\
  && \bigl( \frac {\ikp-\ik}\tau , v \bigr)
  + \bigl( \ka(\nkp)\nabla\ikp , \nabla v \bigr)
  + \kaz (\ikp , v)
  = (\gi , v)
  \label{terzakbis}
  \\
  && \bigl( \frac {\hkp-\hk}\tau , v \bigr)
  + \bigl( \ka(\nkp)\nabla\hkp , \nabla v \bigr)
  + \kaz (\hkp , v)
  = (\gh , v)
  \label{quartakbis}
\Esist
all for every $v\in V$, where we have set
\Bsist
  && \gn := (\alpha-\mu-\oldgianni\ik+\kaz) \nkp 
  \non
  \\
  && \gs := \alpha\nkp - A(\nkp)\skp\ik - A(\nkp)\skp(\hk-\sk) + (\kaz-\mu)\skp
  \non
  \\
  && \gi := \hkp-\skp + (\kaz-1-\mu) \ikp \oldgianni{{} - \nkp \ikp}
  \non
  \\
  && \gh := \alpha\nkp + \skp + (\kaz-1-\mu)\hkp .
  \non
\Esist
We observe that \eqref{basic}, \accorpa{sinfty}{iinfty} and the continuity of $A$ imply that
\Bsist
  && \norma\gn
  + \norma\gs
  + \norma\gi
  + \norma\gh
  \non
  \\
  && \leq c \, \bigl(\norma\gn_\infty
  + \norma\gs_\infty
  + \norma\gi_\infty
  + \norma\gh_\infty \bigr)
  \leq c \,.
  \label{rhsbdd}
\Esist
Now, we test \eqref{primakbis} by $\tau\nkp$ and sum over $k$ from $0$ to $m-1$, 
with an arbitrary positive integer $m\leq N$.
We obtain
\Beq
  \somma k0{m-1} (\nkp-\nk,\nkp)
  + \tau \somma k0{m-1} \bigl (\ka(\nkp)\nabla\nkp,\nabla\nkp) + \kaz \norma\nkp^2 \bigr)
  = \tau \somma k0{m-1} (\gn , \nkp) .
  \non
\Eeq
Since $\ka\geq\kaz$, by owing to \eqref{elementare} and the Schwarz and Young inequality,
we deduce~that
\Bsist
  && \frac 12 \, \norma{\nm}^2
  - \frac 12 \, \norma{n_0}^2
  + \frac 12 \somma k0{m-1} \norma{\nkp-\nk}^2
  + \kaz \tau \somma k0{m-1} \normaV\nkp^2
  \non
  \\
  && \leq \frac \kaz 2 \, \tau \somma k0{m-1} \norma\nkp^2
  + \frac 1 {2\kaz} \, \tau \somma k0{m-1} \norma\gn^2 \,.
  \non
\Esist
Then, it suffices to rearrange and account for the first \oldgianni{condition in} \eqref{stimadatitau} and \eqref{rhsbdd} to~obtain
\Beq
  \max_{k=0,\dots,m} \norma\nk^2
  + \tau \somma k0N \normaV\nk^2
  + \somma k0{N-1} \norma{\nkp-\nk}^2
  \leq c \,.
  \label{stiman}
\Eeq
By treating \accorpa{secondakbis}{quartakbis} in the same way, we also have
\Bsist
  && \max_{k=0,\dots,m} \norma\sk^2
  + \tau \somma k0N \normaV\sk^2
  + \somma k0{N-1} \norma{\skp-\sk}^2
  \leq c
  \label{stimas}
  \\
  && \max_{k=0,\dots,m} \norma\ik^2
  + \tau \somma k0N \normaV\ik^2
  + \somma k0{N-1} \norma{\ikp-\ik}^2
  \leq c
  \label{stimai}
  \\
  && \max_{k=0,\dots,m} \norma\hk^2
  + \tau \somma k0N \normaV\hk^2
  + \somma k0{N-1} \norma{\hkp-\hk}^2
  \leq c \,.
  \label{stimah}
\Esist

At this point, we introduce the interpolating functions. 
Namely, we term
\Beq
  \hn, \hs,\ \hi,\ \hh,\quad
  \on, \os,\ \oi,\ \oh
  \aand
  \un, \us,\ \ui,\ \uh
  \label{interpolating}
\Eeq
the functions $\hat z$, $\overline z$ and $\underline z$ given in Notation~\ref{Interpol}
with $z=(n_0,\dots,n_N)$, $z=(s_0,\dots,s_N)$, $z=(i_0,\dots,i_N)$ and $z=(h_0,\dots,h_N)$, respectively.
Then, we can owe to Proposition~\ref{Propinterp} 
to rewrite the uniform bounds \eqref{basic} and \accorpa{sinfty}{hinfty}, which read
\Bsist
  && \norma\on_{\LQ\infty}
  + \norma\os_{\LQ\infty}
  + \norma\oi_{\LQ\infty}
  + \norma\oh_{\LQ\infty}
  \leq c 
  \label{stimao}
\Esist
and the analogues for $\un$, $\us$, $\ui$ and $\uh$.
Next, the estimates \accorpa{stiman}{stimah} imply
\Bsist
  && \norma\on_{\L2V}
  + \norma\os_{\L2V}
  + \norma\oi_{\L2V}
  + \norma\oh_{\L2V}
  \leq c
  \label{stimeoL2V}
  \\
  && \norma{\on-\un}_{\L2H}
  + \norma{\os-\us}_{\L2H}
  \non
  \\
  && \quad {}
  + \norma{\oi-\ui}_{\L2H}
  + \norma{\oh-\uh}_{\L2H}
  \leq c \tau^{1/2} .
  \label{stimeo-u}
\Esist
\oldgianni{Notice that \eqref{stimeo-u} and \eqref{diffLdue} provide the same estimate for $\on-\hn$.}
Moreover, \oldgianni{the property} \eqref{hzLdueZ} with $Z=V$, 
the second \oldgianni{inequality in} \eqref{stimadatitau} and \eqref{stimao} also yield
\Beq
  \norma\hn_{\L2V}
  + \norma\hs_{\L2V}
  + \norma\hi_{\L2V}
  + \norma\hh_{\L2V}
  \leq c \,.
  \label{stimehL2V}
\Eeq
In order to obtain an estimate for the time derivatives and let $\tau$ tend to zero,
we write equations \accorpa{primak}{quartak} in terms of the interpolating function.
We have
\Bsist
  && (\dt\hn , v)
  + (\ui\on , v)
  + (\ka(\on)\nabla\on , \nabla v)
  = (\alpha-\mu) (\on , v)
  \label{primatau}
  \\
  && (\dt\hs , v)
  + (A(\on)\os\ui , v)
  + (A(\on)\oldgianni\os(\uh-\us) , v)
  + \oldgianni{\mu (\os , v)}
  \qquad
  \non
  \\
  && \quad {}
  + (\ka(\on)\nabla\os , \nabla v)
  = \alpha (\on , v)
  \label{secondatau}
  \\
  && (\dt\hi , v)
  + (1+\mu)(\ui , v)
  + (\ka(\on)\nabla\oi , \nabla v)
  = (\oh-\os , v)
  \label{tezatau}
  \\
  && (\dt\hh , v)
  + (1+\mu)(\oh , v)
  + (\ka(\on)\nabla\oh , \nabla v)
  = (\alpha\on+\os , v)
  \label{quartatau}
\Esist
all of them being satisfied \oldgianni{for every $v\in V$ \aet}.
Moreover, requirements \eqref{segnik} and the initial conditions read
\Bsist
  && \on \geq 0 , \quad
  \oh \geq \os \geq 0
  \aand
  \oi \geq 0
  \quad \aeQ
  \label{segnio}
  \\
  && \hn(0) = \nzt \,, \quad
  \hs(0) = \szt \,, \quad
  \hi(0) = \izt 
  \aand
  \hh(0) = \hzt \,.
  \label{cauchytau}
\Esist
\oldgianni{Finally}, the bounds given by \eqref{basic} and \accorpa{sinfty}{iinfty} imply
\Beq
  \nmin \leq \on \leq \nmax , \quad
  \os \leq \smax , \quad
  \oi \leq \imax
  \aand
  \oh \leq \hmax
  \quad \aeQ .
  \label{bounds}
\Eeq

\step
Second a priori estimate

We take an arbitrary $v\in\L2V$, write \eqref{primatau} at the time $t$ and test it by~$v(t)$.
Then, we integrate over~$(0,T)$ and rearrange.
\oldgianni{Thanks} to the above estimates, we obtain
\Beq
  \intQ \dt\hn \, v
  = - \intQ \ka(\on) \nabla\on \cdot \nabla v
  + \intQ (\alpha-\mu-\ui) \on v
  \leq c \norma v_{\L2V} \,.
  \non
\Eeq
By analogously treating the other equations, we conclude that
\Beq
  \norma{\dt\hn}_{\L2\Vp}
  + \norma{\dt\hs}_{\L2\Vp}
  + \norma{\dt\hi}_{\L2\Vp}
  + \norma{\dt\hh}_{\L2\Vp}
  \leq c \,.
  \label{stimederivate}
\Eeq

Since all this holds under the only restriction on $\tau$ given by~\eqref{hptauz}, 
we are ready to let $\tau$ tend to zero.

\step
Convergence and conclusion

The estimates obtained before and the Aubin-Lions lemma (see, e.g., \cite[Thm.~5.1, p.~58]{Lions})
ensure the existence of a quadruple $\soluz$ such~that
\Bsist
  && \hn \to n \,, \quad
  \hs \to s \,, \quad
  \hi \to i 
  \aand
  \hh \to h 
  \non
  \\
  && \quad \hbox{weakly star in $\H1\Vp\cap\C0H\cap\L2V$}
  \qquad
  \non
  \\
  && \quad \hbox{strongly in $\L2H$ and \aeQ}
  \label{convhat}
\Esist
as $\tau\searrow0$ (more precisely for a suitable sequence $\tau_j\searrow0$).
In particular, the values at~$0$ converge weakly in $H$ to the \lhs s of \eqref{cauchy} 
and the initial conditions \eqref{cauchy} themselves are satisfied on account of~\eqref{convdatitau}.
Moreover, the differences $\on-\un$ and $\on-\hn$ converge to zero strongly in $\L2H$,
and also \aeQ\ without loss of generality,
so that $\on$ and $\un$ tend to $n$ strongly in $\L2H$.
Since the same holds for the other variables, 
all the inequalities \oldgianni{in} \eqref{segni} are satisfied too.
Moreover, \eqref{bounds} yield
\Beq
  \nmin \leq n \leq \nmax , \quad
  s \leq \smax , \quad
  i \leq \imax
  \aand
  h \leq \hmax
  \quad \aeQ .
  \non  
\Eeq
By combining with the point-wise convergence, we infer~that
\Beq
  \hbox{$\on$, $\os$, $\oi$ and $\oh$
    converge to their limits strongly in $\LQ p$ for every $p\in[1,+\infty)$}
  \non
\Eeq
and the same holds for $\un$, $\us$, $\ui$ and $\uh$.
We have some consequences. 
First, the products like $\ui\on$ in \eqref{primatau} converge to the right products
(in~this case the limit is~$in$)
in~$\LQ p$ for $p\in[1,+\infty)$.
Next, $A(\on)$~converges to $A(n)$ in the same topology,
since it is bounded in $\LQ\infty$ and converges to $A(n)$ \aeQ,
all this for $A$ is continuous.
By the same reason, $\ka(\on)$ converges to $\ka(n)$ strongly in~$\LQ p$ for $p<+\infty$.
From this, we \oldgianni{claim} that
\Bsist
  && \ka(\on)\nabla\on \to \ka(n)\nabla n , \quad
  \ka(\on)\nabla\os \to \ka(n)\nabla s , \quad
  \ka(\on)\nabla\oi \to \ka(n)\nabla i 
  \non
  \\
  && \aand
  \ka(\on)\nabla\oh \to \ka(n)\nabla h 
  \quad \hbox{weakly in $(\L2H)^3$}.
  \non
\Esist
\oldgianni{In fact we} prove the first one, only, since the others are analogous.
For the product, we have weak convergence to $\ka(n)\nabla n$ in $(\LQ q)^3$ for $q\in[1,2)$.
On the other hand, $\ka(\on)$ is bounded in $\LQ\infty$ and $\nabla\on$ is bounded in $(\L2H)^3$,
so that the product is bounded in $(\L2H)^3$ and thus has a weak limit in this topology.
Clearly, this weak limit has to be $\ka(n)\nabla n$.

All this allows us to let $\tau$ tend to zero in \accorpa{primatau}{quartatau}.
We consider the first one.
We write its integrated version, namely
\Beq
  \intQ \bigl(
    \dt\hn \, v
    + \ui\on \, v
    + \ka(\on)\nabla\on \cdot \nabla v
  \bigr)
  = (\alpha-\mu) \intQ \on \, v
  \quad \hbox{for every $v\in\L2V$}.
  \non
\Eeq
Then, it is clear that we can let $\tau$ tend to zero and obtain \eqref{intprima}.
As the same argument works for the other equations, the proof is complete.

%%%%%%%%%%%%%%%%%%%%%%%%%%%%%%%%%%%%%%%%%%%%%%%%%%%%%%%%%%%%%%%%%%%%%%%%

\section{Uniqueness}
\label{UNIQUENESS}
\setcounter{equation}{0}

In this section we give \oldgianni{the proofs of Theorems \ref{Uniqueness1} and~\ref{Uniqueness2}.}

\subsection{Proof of Theorem \ref{Uniqueness1}}
\label{UNIQUENESS1}

We assume that $A$ is locally \Lip\ continuous and that $\ka$ is a positive constant.
We pick any two solutions $(n_j,s_j,i_j,h_j)$, $j=1,2$, 
and prove that they are the same.
First, we make some observations.
By recalling that $n_j$ are bounded (as~well as the other components) and that $\inf n_j>0$,
we fix an interval $[n_\star,n^\star]\subset(0,+\infty)$ that contains all the values of $n_1$ and $n_2$.
Then, the restriction of $A$ to $[n_\star,n^\star]$ is \Lip\ continuous.
As for the estimates we are going to prove,
we reinforce the convention given in Notation~\ref{Constants} for the constants
by allowing the values of $c$ to depend on the fixed solutions (through their $L^\infty$ norms), in addition.

We set for brevity $n:=n_1-n_2$, $s:=s_1-s_2$, $i:=i_1-i_2$ and $h:=h_1-h_2$.
We write \eqref{prima} for both solutions, take the difference and test it by~$n$.
Then, we integrate over~$(0,t)$ and adjust.
By also owing to the Young inequality, we~have
\Bsist
  && \frac 12 \, \iO |n(t)|^2
  + \ka \intQt |\nabla n|^2
  = (\alpha-\mu) \intQt |n|^2
  - \intQt \oldgianni{n_1 i n}
  - \intQt i_2 n^2 
  \non
  \\
  && \leq c \intQt \bigl( |n|^2 + |i|^2 \bigr) \,.
  \label{uprima}
\Esist
By proceeding with \eqref{seconda} in the same way, we obtain
\Bsist
  && \frac 12 \, \iO |s(t)|^2
  + \ka \intQt |\nabla s|^2
  + \mu \intQt |s|^2
  \non
  \\
  && =  \alpha \intQt n s
  - \intQt \bigl(
    (A(n_1)-A(n_2)) s_1 i_1
    + A(n_2) s \, i_1
    + A(n_2) s_2 i
  \bigr) s
  \non
  \\
  && \quad{}
  - \intQt \bigl(
    (A(n_1)-A(n_2)) s_1 (s_1-h_1)
    + A(n_2) s (s_1-h_1)
    + A(n_2) s_2 (s-h)
  \bigr) s
  \qquad
  \non
  \\
  && \leq c \intQt \bigl( |n|^2 + |s|^2 + |i|^2 + \oldgianni{|h|^2}\bigr) \,.
  \label{useconda}
\Esist
Analogously, by dealing with \eqref{terza} and \eqref{quarta}, we have
\Bsist
  && \frac 12 \, \iO |i(t)|^2
  + \ka \intQt |\nabla i|^2
  + \frac 12 \, \iO |h(t)|^2
  + \ka \intQt |\nabla h|^2
  \non
  \\
  && \leq c \intQt \bigl( |n|^2 + |s|^2 + |i|^2 + |h|^2 \bigr) \,.
  \label{ualtre}
\Esist
At this point, it suffices to add \accorpa{uprima}{ualtre} to each other
and apply the Gronwall lemma to conclude that $n=s=i=h=0$.

\subsection{Proof of Theorem \ref{Uniqueness2}}
\label{UNIQUENESS2}

It is understood that the assumptions listed in the statement are in force.
Since the component $n$ of every solution belongs to $\Linfty$ 
and is bounded away from zero (see \Regsoluz),
we can assume, without loss of generality, that 
\Beq
  \kaz \leq \ka(y) \leq \kau
  \quad \hbox{for some positive constants $\kaz$ and $\kau$ and every $y>0$}
  \label{kabounds}
\Eeq
and that $A$ is \Lip\ continuous, whenever we fix one or two solutions.
The bounds $\kaz$ and $\kau$ and the \Lip\ constant of $A$ 
depend on the solutions we fix every time, 
and the same holds for the \Lip\ constants of the functions we are going to introduce.
We~set
\Beq
  K(y) := \int_0^y \ka(z) \, dz
  \quad \hbox{for $y\in(0,+\infty)$}
  \label{defKbis}
\Eeq
and observe that, under condition \eqref{kabounds}, it is a bijection from $(0,+\infty)$ onto itself and that
both $K$ and $K^{-1}$ are \Lip\ continuous.
Next, to every solution $\soluz$ to problem \Pbl\ satisfying \Regsoluz\
we associate the function $u$ defined~by
\Beq
  u := K(n) 
  \label{defu}
\Eeq
and observe that the regularity and boundedness properties of~$u$,
but that of the time derivative, are the same as those of~$n$.

Our project is the following:
first we prove that every solution $\soluz$ to problem \Pbl\ satisfying \Regsoluz\
enjoys the regularity properties specified in~\eqref{highreg};
then, we prove uniqueness on account of the regularity already proved.
The main difficulty is a regularity result for the component $n$ of the solution.
This needs some preliminary results.

\Bprop
\label{Regnu}
The component $n$ of every solution $\soluz$ to problem \Pbl\ satisfying \Regsoluz\
and the corresponding function $u$ verify
\Bsist
  && \dt n \in \L2H
  \aand
  \dt u \in \L2H
  \label{regdt}
  \\
  && \div(\ka(n)\nabla n) \in \L2H 
  \aand
  \non
  \\
  && \quad \dt n - \div(\ka(n)\nabla n)
  = (\alpha - \mu - i) n
  \quad \aeQ
  \quad \hbox{with} \quad \hbox{$\dn n=0$ on $\Gamma$}
  \label{bvpn}
  \\
  && \Delta u \in \L2H \aand
  \non
  \\
  && \quad \dt u - \ka(n) \Delta u
  = \ka(n) (\alpha - \mu - i) n
  \quad \aeQ
  \quad \hbox{with} \quad \hbox{$\dn u=0$ on $\Gamma$}.
  \qquad
  \label{bvpu}
\Esist
\Eprop

\Bdim
We fix a solution $\soluz$ 
and suitably extend the component~$n$.
We also introduce an auxiliary function.
We define $\tilde n$ and $f$ on $(-1,T)$ by setting
\Bsist
  && \tilde n(t) := n(t)
  \aand
  f(t) := (\alpha-mu-i(t)) n(t)
  \quad \hbox{if $t\in(0,T)$}
  \non
  \\
  && \tilde n(t) := \nz 
  \aand
  f(t) := - \div (\ka(\nz)\nabla\nz)
  \quad \hbox{if $t\in(-1,0)$}.
  \non
\Esist
Notice that $\div(\ka(\nz)\nabla\nz)$ belongs to~$H$ by \eqref{hpdatareg} 
(for this it would be sufficient a much weaker assumption).
Also notice that $\tilde n$ is a continuous $H$-valued function
so that no jumps for its derivative can occur at $t=0$.
Therefore, $\tilde n$ and $f$ satisfy
\Bsist
  && \tilde n \in H^1(-1,T;\Vp) \cap L^2(-1,T;V) 
  \aand
  f \in L^2(-1,T;H)
  \non
  \\
  && \< \dt \tilde n , v >
  + \iO \ka(\tilde n) \nabla\tilde n \cdot \nabla v
  = \iO f v
  \quad \hbox{for every $v\in V$ a.e.\ in $(-1,T)$}.
  \label{proleqn}
\Esist
However, for simplicity, we write $n$ in place of $\tilde n$ in the sequel.
Now, we fix $\tau>0$ small (namely, $\tau<\min\{1,T\}$). 
Later on, we let $\tau$ tend to zero.
We~set
\Beq
  u(t) := K(n(t))
  \aand
  U(t) := \int_{-\tau}^t u(y) \, dy
  \quad \hbox{for $t\in[-\tau,T]$} 
  \non
\Eeq
by noting that $\nabla u=\ka(n)\nabla n$ and $\dt U=u$.
Then, for $t\in(-\tau,T-\tau)$, we integrate \eqref{proleqn} on $(t,t+\tau)$ with respect to time
and test the resulting equality by $\tau^{-2}(u(t+\tau)-u(t))$.
We obtain for a.a.\ $t\in(-\tau,T-\tau)$
\Bsist
  && \iO \frac{n(t+\tau) - n(t)} \tau \, \frac{u(t+\tau) - u(t)} \tau
  + \iO \nabla \, \frac{U(t+\tau) - U(t)} \tau \cdot \nabla \, \frac{u(t+\tau) - u(t)} \tau
  \non
  \\
  && = \iO \Bigl(  {\textstyle\frac 1\tau \int_t^{t+\tau} f(t') \, dt'} \Bigr) \, \frac{u(t+\tau) - u(t)} \tau \,.
  \non
\Esist
We treat each term of this equality, separately.
Since $n=K^{-1}(u)$ and the derivative of $K^{-1}$ is bounded from below by~$1/\kau$, we~have
\Beq
  \iO \frac{n(t+\tau) - n(t)} \tau \, \frac{u(t+\tau) - u(t)} \tau
  \geq \frac 1\kau \iO \Bigl| \frac{u(t+\tau) - u(t)} \tau \Bigr|^2 .
  \non
\Eeq
As for the second term, we recall that $u=\dt U$, whence
\Beq
  \iO \nabla \, \frac{U(t+\tau) - U(t)} \tau \cdot \nabla \, \frac{u(t+\tau) - u(t)} \tau
  = \frac 12 \, \frac d{dt} \iO \Bigl| \nabla \, \frac{U(t+\tau) - U(t)} \tau \Bigr|^2 .
  \non
\Eeq
Finally, the Young inequality yields
\Bsist
  && \iO \Bigl(  {\textstyle\frac 1\tau \int_t^{t+\tau} f(t') \, dt'} \Bigr) \, \frac{u(t+\tau) - u(t)} \tau
  \non
  \\
  && \leq \frac 1 {2\kau} \iO \Bigl| \frac{u(t+\tau) - u(t)} \tau \Bigr|^2 
  + \frac \kau 2 \iO \Bigl| {\textstyle\frac 1\tau \int_t^{t+\tau} f(t') \, dt'} \Bigr|^2 .
  \non
\Esist
By collecting all this, rearranging and integrating over $(-\tau,T-\tau)$ we obtain
\Bsist 
  && \frac 1\kau \int_{-\tau}^{T-\tau} \iO \Bigl| \frac{u(t+\tau) - u(t)} \tau \Bigr|^2 dt
  + \iO \Bigl| \nabla \, \frac{U(T) - U(T-\tau)} \tau \Bigr|^2
  \non
  \\
  && \leq \iO \Bigl| \nabla \, \frac{U(0) - U(-\tau)} \tau \Bigr|^2
  + \kau \int_{-\tau}^{T-\tau} \iO \Bigl| {\textstyle\frac 1\tau \int_t^{t+\tau} f(t') \, dt'} \Bigr|^2 dt \,.
  \non
\Esist
Now, we observe that the last term can be estimate  uniformly with respect to $\tau$
by the norm of $f$ and $\nz$ in $\L2H$ and~$\Wx{2,\infty}$, respectively,
and that $U(0)-U(-\tau)=\tau K(\nz)$.
Therefore, we conclude that $\dt u\in\L2H$.
As $n=K^{-1}(u)$ and $K^{-1}$ is \Lip\ continuous, we infer that $\dt n\in\L2H$,
so that \eqref{regdt} is proved.
This conclusions also imply that $\dt u=\ka(n)\dt n$ \aeQ.
Next, we come back to \eqref{proleqn} written on $(0,T)$ and recall that $f\in\LQ\infty$.
Thus, \eqref{bvpn} follows.
Finally, we can rewrite \eqref{proleqn} on $(0,T)$ this~way
\Beq
  \iO \frac 1{\ka(n)} \, \dt u \, v 
  + \iO \nabla u \cdot \nabla v
  = \iO f v
  \quad \hbox{for every $v\in V$ a.e.\ in $(0,T)$}.
  \non
\Eeq
This implies that $\Delta u\in\L2H$ and~that
\Beq
  \frac 1{\ka(n)} \, \dt u - \Delta u = f
  \quad \aeQ
  \aand
  \dn u = 0
  \quad \hbox{on $\Gamma$}.
  \non
\Eeq
Then \eqref{bvpu} follows.
\Edim

\Blem
\label{Unic}
Let $a\in\CQ0$ be strictly positive and let $z$ satisfy
\Bsist
  && z \in \H1H \cap \L2W
  \non
  \\
  && a \dt z - \Delta z = 0 \quad \aeQ
  \aand 
  z(0) = 0 .  
\Esist
Then $z=0$.
\Elem

\Bdim
By multiplying the equation by $\dt z$ and integrating over~$Q_t$, we obtain
\Beq
  \intQt a |\dt z|^2
  + \frac 12\iO |\nabla z(t)|^2
  = 0
  \quad \hbox{for every $t\in[0,T]$}.
  \non
\Eeq
Since $a$ is strictly positive, we deduce that 
both $\dt z$ and $\nabla z$ vanish \aeQ, whence $z$ is a constant.
As $z(0)=0$, we conclude that $z=0$.
\Edim

\Bprop
\label{Highreg}
Every solution to problem \Pbl\ satisfying \Regsoluz\
also enjoys the regularity specified by~\eqref{highreg}.
\Eprop

\Bdim
Fix a solution $\soluz$.
First, we consider the component~$n$ and recall that
it satisfies \accorpa{regdt}{bvpn}, where the \rhs\ of the equation belongs to $\LQ\infty$.
Since \eqref{hpdatareg} implies $\nz\in\Cx0$, 
by applying \cite[Thm.~1.3 and Rem.~1.1 of Chpt.~III]{DiB},
we deduce that $n\in\CQ0$, whence also $\ka(n)\in\CQ0$.
Now, for $p\in[2,+\infty)$, we consider the problem of finding $w$ satisfying
\Bsist
  && w \in \H1H \cap \L2W
  \label{regw}
  \\
  && \dt w - \ka(n) \Delta w
  = \ka(n) (\alpha - \mu - i) n
  \quad \aeQ
  \aand
  w(0) = K(\nz) \,.
  \qquad
  \label{eqw}
\Esist
Thanks to Proposition~\ref{Regnu}, $u$~is a solution.
On the other hand, $\ka(n)$ is continuous as just observed,
the \rhs\ of the equation is bounded and $K(\nz)$~belongs to~$\Wx{2,\infty}$ and satisfies $\dn K(\nz)=0$ on~$\Gamma$
(recall that $K$ is of class~$C^2$ with bounded derivatives by \eqref{hpnonlin} and \eqref{kabounds}).
Therefore, we can apply \cite[Thm.~2.1]{DHP} and ensure the existence 
of a unique $w\in\W{1,p}{\Lx p}\cap\L p{\Wx{2,p}}$ 
that satisfies \eqref{eqw} and the homogeneous Neumann boundary condition.
Since $p\geq2$, \eqref{regw} holds as well 
and the assumptions of Lemma~\ref{Unic} with $a=1/\ka(n)$ are fulfilled by $w-u$.
We conclude that $w=u$, so~that
\Beq
  u \in \W{1,p}{\Lx p} \cap \L p{\Wx{2,p}}.
  \non
\Eeq
Since $p$ is arbitrary in $[2,+\infty)$, $\Omega$ is bounded, 
$n=K^{-1}(u)$ and $K^{-1}$ is of class~$C^2$ 
with bounded derivatives, we deduce that \eqref{highreg} holds for~$n$.

Now, we consider the components $s$, $i$ and~$h$.
We observe that each of equations \accorpa{seconda}{quarta} 
and the corresponding boundary and initial conditions have the form
\Beq
  \dt w - \div(\ka(n)\nabla w) = f
  \quad \hbox{with $\dn w=0$}
  \aand
  w(0) = w_0
  \label{eqdiv}
\Eeq
where $f\in\LQ\infty$ and $w_0\in\Wx{2,\infty}$ with $\dn w_0=0$.
On account of the already proved regularity of~$n$,
we can rewrite the equation in \eqref{eqdiv} in the non-divergence form, namely
\Beq
  \dt w - \bigl( \ka(n) \Delta w + \ka'(n) \nabla n \cdot \nabla w) = f 
  \label{eqnondiv}
\Eeq 
and apply \cite[Thm.~2.1]{DHP} once more with $p\geq2$
(using the $4p$-summability of $\nabla n$).
Then the homogeneous Neumann problem for \eqref{eqnondiv} 
complemented with the initial condition $w(0)=w_0$
has a unique solution $w\in\W{1,p}{\Lx p}\cap\L p{\Wx{2,p}}$.
Clearly, such a $w$ belongs to $\H1H\cap\L2W$ and satisfies \eqref{eqdiv}.
Since \eqref{eqdiv} has a unique solution in $\H1H\cap\L2W$,
we deduce that $w\in\W{1,p}{\Lx p}\cap\L p{\Wx{2,p}}$
whenever $w$ belongs to $\H1H\cap\L2W$ and satisfies \eqref{eqdiv}.
This is the case for $s$, $i$ and~$h$.
Since $p\geq2$ is arbitrary and $\Omega$ is bounded, \eqref{highreg} is completely proved.
\Edim

\step
Conclusion of the proof of Theorem~\ref{Uniqueness2}

We pick any two solutions $(n_j,s_j,i_j,h_j)$, $j=1,2$, 
and prove that they are the same.
Thanks to Proposition~\ref{Highreg}, they satisfy \eqref{highreg}.
In particular, the gradients of all the components belong to $\L4\Linfty$
since $\Wx{2,4}\subset\Wx{1,\infty}$.
For simplicity, we use the rule of Notation~\ref{Constants} concerning the constants:
the symbol $c$ stands for possible difference constants 
that depend on the structure, the data, some norms of the solutions we have fixed
and the constants appearing in \eqref{kabounds} 
(which are chosen after fixing the solutions) 
and the consequent \Lip\ constants of the nonlinearities.
As in the proof of Theorem~\ref{Uniqueness1}, we set for brevity 
$n:=n_1-n_2$, $s:=s_1-s_2$, $i:=i_1-i_2$ and $h:=h_1-h_2$,
write equations \accorpa{prima}{quarta} for both solutions and test the differences by
$n$, $s$, $i$ and~$h$, respectively.
As for the first equation, we~have
\Bsist
  && \frac 12 \, \iO |n(t)|^2
  + \intQt \ka(n_1) |\nabla n|^2
  \non
  \\
  && = \intQt (\ka(n_1)-\ka(n_2)) \nabla n_2 \cdot \nabla n
  + (\alpha-\mu) \intQt |n|^2
  - \intQt n_1 i n
  - \intQt i_2 n^2 \,.
  \non
\Esist
While we can simply use the inequality $\ka(n_1)\geq\kaz$ on the \lhs,
the true novelty is the first term on the \rhs.
We treat it by owing to the \Holder\ and Young inequalities this~way
\Bsist
  && \intQt (\ka(n_1)-\ka(n_2)) \nabla n_2 \cdot \nabla n
  \leq c \iot \norma{n(t')} \, \norma{\nabla n_2(t')}_\infty \, \norma{\nabla n(t')} \, dt'
  \non
  \\
  && \leq \kaz \intQt |\nabla n|^2
  + c \iot \norma{\nabla n_2(t')}_\infty^2 \, \norma{n(t')}^2 \, dt' .
  \non 
\Esist
The other terms on the \rhs\ can be dealt with as in the above proof.
By combining and rearranging, we easily deduce that
\Beq
  \iO |n(t)|^2
  \leq c \iot \bigl( \norma{\nabla n_2(t')}_\infty^2 +1 \bigr) \bigl( \norma{n(t')}^2 + \norma{i(t')}^2 \bigr) \, dt'\,.
  \label{okprima}
\Eeq
The next equation \eqref{seconda} can be treated in a similar way.
Arguing as we did to derive \eqref{okprima} and~\eqref{useconda},
we~have
\Bsist
  && \frac 12 \, \iO |s(t)|^2
  + \intQt \ka(n_1) |\nabla s|^2
  + \mu \intQt |s|^2
  \non
  \\
  && = \intQt (\ka(n_1)-\ka(n_2) \nabla s_2 \cdot \nabla s
  \non
  \\
  && \quad {}
  + \alpha \intQt n s
  - \intQt \bigl(
    (A(n_1)-A(n_2)) s_1 i_1
    + A(n_2) s \, i_1
    + A(n_2) s_2 i
  \bigr) s
  \non
  \\
  && \quad{}
  - \intQt \bigl(
    (A(n_1)-A(n_2)) s_1 (s_1-h_1)
    + A(n_2) s (s_1-h_1)
    + A(n_2) s_2 (s-h)
  \bigr) s
  \qquad
  \non
  \\
  && \leq \kaz \intQt |\nabla s|^2
  + c \iot \bigl( \norma{\nabla s_2(t')}_\infty^2 +1 \bigr)
  \bigl(
    \norma{n(t')}^2 + \norma{s(t')}^2 + \norma{i(t')}^2 + \norma{h(t')}^2
  \bigr) \, dt' \,.
  \non
\Esist
Moreover, we can use the inequality $\ka(n_1)\geq\kaz$ on the \lhs\ and rearrange also in this case.
By treating equations \accorpa{terza}{quarta} in the same way and summing up,
we conclude~that
\Bsist
  && \iO \bigl( |n(t)|^2 + |s(t)|^2 + |i(t)|^2 + |h(t)|^2 \bigr)
  \non
  \\
  && \leq \iot \psi(t')
  \bigl(
    \norma{n(t')}^2 + \norma{s(t')}^2 + \norma{i(t')}^2 + \norma{h(t')}^2
  \bigr) \, dt'
  \quad \hbox{for every $t\in[0,T]$}
  \non
\Esist
with a function $\psi\in L^1(0,T)$.
Thus, the Gronwall lemma yields $n=s=i=h=0$,
and the proof is complete.

%%%%%%%%%%%%%%%%%%%%%%%%%%%%%%%%%%%%%%%%%%%%%%%%%%%%%%%%%%%%%%%%%%%%%%%%

\section*{Acknowledgments}
\pier{This research was supported by the Italian Ministry of Education, 
University and Research~(MIUR): Dipartimenti di Eccellenza Program (2018--2022) 
-- Dept.~of Mathematics ``F.~Casorati'', University of Pavia. 
In addition, {PC and ER gratefully mention} some other support 
from the MIUR-PRIN Grant 2020F3NCPX 
``Mathematics for industry 4.0 (Math4I4)'' and the GNAMPA (Gruppo Nazionale per l'Analisi Matematica, 
la Probabilit\`a e le loro Applicazioni) of INdAM (Isti\-tuto 
Nazionale di Alta Matematica), while AR gratefully acknowledges the partial support of the MIUR-PRIN project XFAST-SIMS (no. 20173C478N).}

%%%%%%%%%%%%%%%%%%%%%%%%%%%%%%%%%
%% bibliography
%%%%%%%%%%%%%%%%%%%%%%%%%%%%%%%%%

\vspace{3truemm}

\Begin{thebibliography}{10}

\bibitem{Albi} 
\pcol{G. Albi, L. Pareschi, M. Zanella,
Control with uncertain data of socially structured compartmental epidemic models,
{\it J. Math. Biol.}  {\bf 82}  (2021) Paper No. 63, 41 pp.}

\bibitem{Barbu}
V. Barbu,
\gianni{``Nonlinear Differential Equations of Monotone Types in Banach Spaces'',
Springer, 
New York,
2010.}

\bibitem{Bellomo}
\ale{N. Bellomo, R. Bingham, M.A.J. Chaplain, G. Dosi, G. Forni, D.A. Knopoff, J. Lowengrub, R. Twarock, M.E. Virgillito,
A multiscale model of virus pandemic: heterogeneous interactive entities in a globally connected world,
{\it \pcol{Math. Models Methods Appl. Sci.}} {\bf 30} (2020) 1591-1651.}

\bibitem{Bellomo2}
\ale{N. Bellomo, K.J. Painter, Y. Tao, M. Winkler,
Occurrence vs. absence of taxis-driven instabilities in a May-Nowak model for virus infection,
{\it SIAM J. Appl. Math.} \pcol{{\bf 79} (2019) 1990-2010.}}

\bibitem{Bellomo3}
\ale{N. Bellomo, N. Outada, J. Soler, Y. Tao, M. Winkler,
Chemotaxis and Cross-diffusion Models in Complex Environments: 
Models and Analytic Problems Towards a Multiscale Vision,
{\it \pcol{Math. Models Methods Appl. Sci.}}, to appear.}

\bibitem{Bellomo4}
\ale{N. Bellomo, F. Brezzi, M.A.J. Chaplain,
Special Issue on ``Mathematics Towards COVID19 and Pandemic'',
{\it \pcol{Math. Models Methods Appl. Sci.}} {\bf 31} (2021) Issue 12.}

\bibitem{Berestycki}
\ale{H. Berestycki, J.-M. Roquejoffre, L. Rossi,
Propagation of epidemics along lines with fast diffusion,
{\it Bull. Math. Biol.} {\bf 83} (2021) \pcol{Paper No. 2, 34 pp.}}

\bibitem{Bertaglia}
\ale{G. Bertaglia, L. Pareschi,
Hyperbolic compartmental models for epidemic spread on networks with uncertain data: application to the emergence of COVID-19 in Italy,
{\it \pcol{Math. Models Methods Appl. Sci.}} {\bf 31} (2021) 2495-2531.}

\bibitem{Calleri} 
\pcol{F. Calleri, G. Nastasi, V. Romano, 
Continuous-time stochastic processes for the spread of COVID-19 disease simulated 
via a Monte Carlo approach and comparison with deterministic models,
{\it J. Math. Biol.}  {\bf 83}  (2021) Paper No. 34, 26 pp.}

\bibitem{DHP}
R. Denk, M. Hieber, J. Pr\"uss,
Optimal $L^p$-$L^q$-estimates for parabolic boundary value problems with inhomogeneous data,
{\it Math. Z.} {\bf 257} (2007) 193-224.

\bibitem{DiB}
E. DiBenedetto,
``Degenerate Parabolic Equations'',
Springer-Verlag, 
New York, 
1993.

\bibitem{Gatto}
\ale{M. Gatto, E. Bertuzzo, L. Mari, S. Miccoli, L. Carraro, R. Casagrandi, A. Rinaldo,
Spread and dynamics of the COVID-19 epidemic in Italy: effects of emergency containment measures,
{\it Proc. Nat. Acad. Sci.} {\bf 117} (2020) 10484-10491.}

\bibitem{Giordano}
\ale{G. Giordano, F. Blanchini, R. Bruno, P. Colaneri, A. Di Filippo, A. Di Matteo, M. Colaneri,
Modelling the COVID-19 epidemic and implementation of population-wide interventions in Italy,
{\it Nat. Med.} {\bf 26} (2020) 855-860.}

\bibitem{Grave}
\ale{M. Grave, A.L.G.A. Coutinho,
Adaptive mesh refinement and coarsening for diffusion-reaction epidemiological models,
{\it Comput. Mech.} {\bf 67} (2021) 1177-1199.}
%https://doi.org/10.1007/s00466-020-01888-0

\bibitem{Grave2}
\ale{M. Grave, A. Viguerie, G.F. Barros, A. Reali, A.L.G.A. Coutinho,
Assessing the spatio-temporal spread of COVID-19 via compartmental models with diffusion in Italy, USA, and Brazil,
{\it Arch. Comput. Mech. Eng.} {\bf 28} (2021) 4205-4223.}

\bibitem{Guglielmi}
\pcol{N. Guglielmi, E. Iacomini, A. Viguerie,
Delay differential equations for the spatially resolved simulation of epidemics with specific application to COVID-19
{\it Math. Methods Appl. Sci}, to appear, DOI:~10.1002/mma.8068}

\bibitem{Jha}
\ale{P.K. Jha, L. Cao, J.T. Oden, 
Bayesian-based predictions of COVID-19 evolution in Texas using multispecies mixture-theoretic continuum models,
{\it Comput. Mech.} {\bf 66} (2020) 1055-1068.}

\bibitem{Linka}
\ale{K. Linka, P. Rahman, A. Goriely, E. Kuhl,
Is it safe to lift COVID-19 travel bans? The Newfoundland story,
{\it Comput. Mech.} {\bf 66} (2020) 1081-1092.}

\bibitem{Lions}
J.-L.~Lions,
``Quelques m\'ethodes de r\'esolution des probl\`emes
aux limites non lin\'eaires'',
Dunod; Gauthier-Villars, Paris, 1969.

\bibitem{Parolini}
\pcol{N. Parolini, L. Ded\`e, P.F. Antonietti, G.  Ardenghi, A. Manzoni, E. Miglio,  A. Pugliese, M.  Verani, A. Quarteroni,  
SUIHTER: a new mathematical model for COVID-19. Application to the analysis of the second epidemic outbreak in Italy,
{\it Proc. A.}  {\bf 477}  (2021), Paper No. 20210027, 21 pp.}

\bibitem{Vig1}
A. Viguerie, G. Lorenzo, F. Auricchio, D. Baroli,
T.J.R. Hughes, A. Patton, A. Reali,
Th.E. Yankeelov, A. Veneziani,
Simulating the spread of COVID-19 via a spatially-resolved
susceptible-exposed-infected-recovered-deceased (SEIRD) model with heterogeneous diffusion,
\pier{{\it Appl. Math. Lett.} {\bf 111} (2021) Paper No. 106617, 9 pp.}

\bibitem{Vig2}
A. Viguerie, A. Veneziani, G. Lorenzo, D. Baroli, N. Aretz-Nellesen,
A. Patton, T.E. Yankeelov, A. Reali, T.J.R. Hughes, F. Auricchio,
\pier{Diffusion-reaction compartmental models formulated in a continuum mechanics framework: application to COVID-19, mathematical analysis, and numerical study,
{\it Comput. Mech.} {\bf 66} (2020) 1131-1152.}
%https://doi.org/10.1007/s00466-020-01888-0

\bibitem{Wang}
\ale{Z. Wang, X. Zhang, G.H. Teichert, M. Carrasco-Teja, K. Garikipati,
System inference for the spatio-temporal evolution of infectious diseases: Michigan in the time of COVID-19,
{\it Comput. Mech.} {\bf 66} (2020) 1153-1176.}

\bibitem{Winkler}
\ale{M. Winkler,
Boundedness in a chemotaxis-May-Nowak model for virus dynamics with mildly saturated chemotactic sensitivity
{\it Acta Appl. Math.} \pcol{{\bf 163} (2019) 1-17.}}

\bibitem{Zohdi}
\ale{T.I. Zohdi,
An agent-based computational framework for simulation of global pandemic and social response on planet X,
{\it Comput. Mech.} {\bf 66} (2020) 1195-1209.}

\End{thebibliography}

\End{document}

%%%%%%%%%%%%%%%%%%%%%%%%%%%%%%%%%%%%%%%%%%%%%

\bibitem{Brezis}
H. Brezis,
``Op\'erateurs maximaux monotones et semi-groupes de contractions
dans les espaces de Hilbert'',
North-Holland Math. Stud.
{\bf 5},
North-Holland,
Amsterdam,
1973.

\bibitem{Simon}
J. Simon,
{Compact sets in the space $L^p(0,T; B)$},
{\it Ann. Mat. Pura Appl.~(4)\/} 
{\bf 146} (1987) 65-96.

%%%%%%%%%%%%%%%%%%%%%%%%%%%%%%%%%%%%%%%%%%%%%%

\Blem
\label{Unicn}
Assume that $\soluz$ is a solution to problem \Pbl\ satisfying \Regsoluz.
Set $f:=(\alpha-\mu-i)n$
and consider the problem of finding $w$ such~that
\Bsist
  && w \in \H1\Vp \cap \L2V \cap \LQ\infty
  \quad \hbox{with} \quad
  \inf w > 0 
  \label{regw}
  \\
  && \< \dt w , v >
  + \iO \ka(w) \nabla w \cdot \nabla v
  = \iO f v
  \quad \hbox{for every $v\in V$ \aet}.
  \label{eqw}
\Esist
Then, $n$ is its unique solution.
\Elem

\Bdim
Since \eqref{eqw} with $w=n$ coincides with \eqref{prima},
$w=n$ is a solution.
Now, let $w$ be any solution: we prove that $w=n$.
We recall \eqref{regw}, owe to \eqref{kabounds}, whose bounds now also depend on~$w$,
and term $L$ the \Lip\ constant of~$\ka$.
Let us introduce $\signeps,\modeps\cdot:\erre\to\erre$ by setting
\Bsist
  && \signeps(y) := y/\eps 
  \quad \hbox{if $|y|\leq\eps$}
  \aand
  \signeps(y) := y/|y|
  \quad \hbox{otherwise}
  \non
  \\
  && \modeps y := \int_0^y \signeps(y') \, dy'
  \quad \hbox{for every $y\in\erre$}.
  \non
\Esist
We set $z:=w-n$ and choose $v=\signeps(z)$ in the difference of \eqref{eqw} 
written for $w$ and~$n$.
By integrating over $(0,t)$, we have
\Beq
  \iO \modeps{z(t)}
  + \intQt \ka(w) \signeps'(z) |\nabla z|^2
  = \intQt (\ka(n)-\ka(w)) \nabla n \cdot \nabla z \, \signeps'(z)
  \non
\Eeq
We estimate the second integral on the \lhs\ from blow by owing to~\eqref{kabounds}
and observe that $|y|(\signeps'(y))^{1/2}\leq\eps^{1/2}$ for every $y\in\erre$.
Hence, we~have
\Bsist
  && \iO \modeps{z(t)}
  + \kaz \intQt \signeps'(z) |\nabla z|^2
  \non
  \\
  && \leq L \intQt |z| (\signeps'(z))^{1/2} |\nabla n| \, |\nabla z| (\signeps'(z))^{1/2}
  \non
  \\
  && \leq L\eps^{1/2} \intQt |\nabla n| \, |\nabla z| (\signeps'(z))^{1/2}
  \non
  \\
  && \leq \frac {L\eps^{1/2}} 2 \intQt \signeps'(z) |\nabla z|^2
  + \frac {L\eps^{1/2}} 2 \intQ |\nabla n|^2.
  \non
\Esist
We deduce that
\Beq
  \iO \modeps{z(t)}
  + (\ka_0 - L\eps^{1/2}/2) \intQt \signeps'(z) |\nabla z|^2
  \leq \frac {L\eps^{1/2}} 2 \intQ |\nabla z|^2.
  \non
\Eeq
Since the second term on the \lhs\ is nonnegative if $\eps<(2\kaz/L)^2$,
by letting $\eps$ tend to zero we obtain $z=0$.
\Edim

We take any $p>1$ and $m>0$,
set $v:=\min\graffe{u,m}$ and test \eqref{pblG} by $v^{p-1}$.
With $\delta:=\bz^{1/p'}$ we obtain
\Bsist
  && \bz \norma v_p^p
  \leq \iO a \nabla u \cdot \nabla(v^{p-1})
  + \iO b u v^{p-1}
  = \iO f v^{p-1}
  \non
  \\
  && \leq \norma f_p \norma{v^{p-1}}_{p'}
  = \delta^{-1} \norma f_p \, \delta \norma v_p^{p-1} \,.
  \non
\Esist
By applying the Young inequality, we deduce that
\Beq
  \bz \norma v_p^p
  \leq \frac 1p \, \delta^{-p} \norma f_p^p
  + \frac 1{p'} \, \bz \norma v_p^p \,,
  \quad \hbox{i.e.,} \quad
  \bz \norma v_p^p
  \leq \delta^{-p} \norma f_p^p \,.
  \non
\Eeq
Hence, we have
\Beq
  \bz^{1/p} \norma v_p
  \leq \delta^{-1} \norma f_p
  = \bz^{-1/p'} \norma f_p \,.
  \non
\Eeq
By letting $p$ tend to infinity we infer that $v$ is bounded and that
\Beq
  \norma v_\infty
  \leq \bz^{-1} \norma f_\infty \,,
  \quad \hbox{i.e.,} \quad
  \min\graffe{u,m}
  \leq \bz^{-1} \norma f_\infty 
  \quad \aeO
\Eeq
and we conclude by letting $m$ tend to infinity.